\documentclass[12pt]{article}
\usepackage{amssymb,amsmath}
\usepackage{amsbsy}
\usepackage{comment}
\usepackage{enumerate}
\everymath{\displaystyle}

\usepackage{enumerate}

\usepackage{hyperref}
\hypersetup{
    colorlinks=true,
    linkcolor=blue,
    filecolor=magenta,
    urlcolor=cyan,
    pdftitle={Overleaf Example},
    pdfpagemode=FullScreen,
    }
\everymath{\displaystyle}
\usepackage{xcolor}

\newtheorem{theorem}{Theorem}[section]

\newtheorem{lemma}[theorem]{Lemma}

\newtheorem{definition}[theorem]{Definition}

\newtheorem{remark}[theorem]{Remark}

\newcommand{\cqd}{\hfill $\square$}

\begin{document}

\title{Blow-up and global mild solutions for a Hardy-H\'enon parabolic equation on the Heisenberg group
	
}
\author{Ricardo Castillo\thanks{Facultad de Ciencias, Departamento de Matem\'atica, Universidade del B\'io-B\'io, Concepci\'on, Chile E-mail: {\tt rcastillo@ubiobio.cl}. Supported by ANID-FONDECYT - 11220152}, Ricardo Freire\thanks{Departamento de Matem\'atica, Universidade Estadual do Sudoeste da Bahia, Vitória da Conquista-BA, Brazil. E-mail: {\tt ricardo.silva@uesb.edu.br}. }, and Miguel Loayza\thanks{Departamento de Matem\'atica, Universidade Federal de Pernambuco, Recife-PE, Brazil. E-mail: {\tt miguel.loayza@ufpe.br}. Partially supported by MATH-AMSUD project 21-MATH-03 (CTMicrAAPDEs), CAPES-PRINT 88887.311962/2018-00 (Brazil)}}
\date{}

\maketitle

\begin{abstract}
We are concerned with the existence of global and blow-up solutions for the nonlinear parabolic problem described by the Hardy-Hénon equation $u_t - \Delta_{\mathbb{H}} u = |\cdot|_{\mathbb{H}}^{\gamma} u^p \mbox{ in } \mathbb{H}^N \times (0,T),$ 
where $\mathbb{H}^N$ is the $N$-dimensional Heisenberg group, and the singular term $|\cdot|_{\mathbb{H}}^{\gamma}$ is given by the Korányi norm. Our study focuses on nonnegative solutions. We establish that for $\gamma\geq 0$, the Fujita critical exponent is $p_c = 1+ (2+\gamma)/Q$, where $Q=2N+2$ is the homogeneous dimension of $\mathbb{H}^N$. For $\gamma<0$, the solutions blow up for $1<p<1+ (2+\gamma)/Q$, while global solutions exist for $p>1+ (2+\gamma)/(Q + \gamma )$. In particular,  our results coincide with the results found by Georgiev and Palmieri in \cite{PALMIERI} for $\gamma=0$.\\
~~\\
~~\\	
MSC[2020]{Primary 35A01, 35B44, 35K55, 35K08, 35R03} 
	
\end{abstract}

\section{Introduction}
Let $T>0$, $p>1$ and $\gamma>-2$. We consider the following semilinear parabolic problem
\begin{equation}\label{Eqgeral-1}
	\left\{ 
	\begin{array}{rll}
		u_t -  	\Delta_{\mathbb{H}} u &=  |\cdot|^{\gamma}_{\mathbb{H}} u^p& \mbox{ in } \mathbb{H}^N \times (0,T), \\	
		{u}(0) &=  {u}_{0}& \mbox{ in } \mathbb{H}^N, \\
	\end{array}
	\right.
\end{equation}
where $\mathbb{H}^N$ is the $N$-dimensional Heisenberg group, and $|\cdot|_{\mathbb{H}}$ is the Korányi or an equivalent symmetric homogeneous norm. When $\gamma < 0$, equation \eqref{Eqgeral-1} is known as the Hardy parabolic equation, while for $\gamma > 0$, it is called the Hénon parabolic equation. 

In the Euclidean frame $\mathbb{R}^N$,  the elliptic version of \eqref{Eqgeral-1} is $-\Delta u = |\cdot|^{\gamma} u^p$ which was  introduced by Hénon as a model to study rotating stellar systems (see \cite{Henon}), and has been extensively studied in mathematics, particularly in nonlinear analysis and variational methods (see \cite{Ghou}). In addition, numerous researchers have spent decades determining the so-called Fujita critical exponent  $p_c$ for problems related to \eqref{Eqgeral-1}. The exponent $ p_c $ is critical in the following sense: if $ 1 < p \leq p_c$ , then problem \eqref{Eqgeral-1} does not admit a global solution (blow-up solution), while for $ p > p_c $, global solutions exist (see \cite{FU} and \cite{QS}). It was determined that $p_c = 1 + (2 + \gamma)/N$ is the Fujita critical exponent for the problem \eqref{Eqgeral-1} (when $\mathbb{H}^N$ is replaced by $\mathbb{R}^N$), see \cite{QI}, \cite{SL}, \cite{CHI}, \cite{PINSKY1} (non singular case), and the references therein. 


When $\gamma=0$, problem \eqref{Eqgeral-1} was studied, within the Carnot groups framework, by Andrea Pascucci in \cite{Pascucci1} who showed that the Fujita critical exponent is $p_c = 1+ 2/Q$, where $Q$ represents the homogeneous dimension of the Carnot group. More recently, within the context of Heisenberg groups, Georgiev and Palmieri also derived the same exponent with $Q=2N+2$ (see \cite{PALMIERI}). They used the test function method to prove the nonexistence of global solutions for $p \leq 1 + 2/Q$ and established the global existence for small initial data in a specific class of weighted $L^{\infty}(\mathbb{H}^N)$ spaces when $p > 1 + 2/Q$. Several authors have addressed the existence and nonexistence of solutions relative to the homogeneous dimension of the Heisenberg group; see for instance, \cite{GOLDSTEIN}, \cite{Goldstein1}, \cite{AHMAD}, \cite{PALMIERI}, \cite{GEORGIEV2}, \cite{RUZHANSKY}.
Consequently, the critical exponent for the problem \eqref{Eqgeral-1} is expected to depend on the homogeneous dimension $Q=2N+2$. Pohozaev and Véron confirmed this dependency in their work, as evidenced in \cite{MITIDIERI}, where they obtained nonexistence results for elliptic, parabolic, and hyperbolic equations within the Heisenberg group, including the equation \eqref{Eqgeral-1}. They showed that problem \eqref{Eqgeral-1} does not admit weak non-trivial solutions in weighted space $L^p_{\text{loc}}(\mathbb{R}^{2N+1} \times (0,\infty), |\varsigma|_{\mathbb{H}}^{\gamma} d\varsigma dt)$ when $p < 1 + (2 + \gamma)/Q$, see \cite[Theorem 3.1]{MITIDIERI}.


Motivated by the aforementioned findings, our primary objective is to determine conditions for the global and nonglobal existence of mild solutions to the problem (\ref{Eqgeral-1}). When $\gamma \geq 0$, we prove that $p_c = 1 + (2 + \gamma)/Q$ is the Fujita critical exponent for the problem \eqref{Eqgeral-1}. For $\gamma < 0$, we show that there are no global nontrivial and nonnegative solutions to the problem \eqref{Eqgeral-1} if $1 < p \leq  1 + (2 + \gamma)/Q$ (blow-up solution), and that there is a global solution if $p > 1 + (2 + \gamma)/(Q+\gamma)$. Unlike previous works, to achieve our results, we combined new semigroup estimates (see, for example, Lemma \ref{UpperSemigroup}, Lemma \ref{Henon1}, and inequality (\ref{eqfc3-A1})) with necessary adaptations to the methodology proposed by Weissler and Fujita (in the Euclidean framework, see \cite{Weissler}, \cite{Weissler1}, \cite{Weissler2}, and \cite{FU1}). For example, to study the blow-up, we transform the mild formulation of the solution into a suitable differential inequality (see Lemma \ref{lemma-tipo-fuj}). Several recent research papers have successfully employed this approach in the Euclidean context. For instance, see \cite[Lemma 4.1]{Fujishima1}, \cite[Lemma 2.3]{Fujishima2}, and \cite[Page 444]{CastilloHardy}.



Before presenting the results obtained, we introduce the notion of the solution considered in this work.  For $\gamma\geq 0$, we need to consider the weighted space $$L^{\infty}_{\varphi}(\mathbb{H}^N)= \left \{ \phi \in L^{\infty}(\mathbb{H}^N); \; \| \varphi \cdot \phi \|_{L^{\infty} ( \mathbb{H}^N )}  < +\infty \right \},$$
where $\varphi(\eta) = (1+|\eta|_{\mathbb{H}})^{\gamma/(p-1)}$. This space was used by Wang in \cite{WANG} in the Euclidean context. When $\gamma<0$, we will use the Banach space $L^\infty (\mathbb{H}^N)$.


\begin{definition}\label{mild-solution} Suppose $X$ is either $L^\infty(\mathbb{H}^N)$ or $ L^{\infty}_{\varphi}(\mathbb{H}^N)$ and $u_0 \in X$. A mild solution (or simply a solution) of the problem \eqref{Eqgeral-1} is a function $u \in L^{\infty}([0,T], X)$ that verifies the integral equation
	\begin{equation}\label{defmildsol}
		u(t)= e^{ t\Delta_{\mathbb{H}} }u_0 + \int_{0}^{t} e^{(t-s)\Delta_{\mathbb{H}}} \,  | \cdot |_{\mathbb{H}}^{\gamma} \, [u(s)]^p ds,
	\end{equation}
for a.e. $t\in [0,T]$. This solution can be extended to a maximal interval $[0, T_{\max})$. When $T_{\max}  =\infty$, we say that $u$ is a global solution of \eqref{Eqgeral-1}. Otherwise, we say that the solution blows up in finite time; see Section \ref{Exist}.

\end{definition}

The following result shows that for $\gamma \geq  0$, the Fujita critical exponent for problem \eqref{Eqgeral-1} is given by $p_c = 1 + (2 + \gamma)/Q$.

\begin{theorem}\label{global-e-blowup-gamapositivo} Suppose $p>1$ and $\gamma\geq 0$.
    \begin{enumerate}[(i)]
        \item Let $0\leq  \gamma <Q(p-1)$. If
$1<p \leq p_c = 1 + {(2+\gamma)}/{Q}$, then all nontrivial nonnegative solutions of problem \eqref{Eqgeral-1} blow up in finite time.
 
 \item If $p_c= 1 + {(2+\gamma)}/{Q} <p$,
then there are nontrivial and nonnegative initial data in $L^{\infty}_{\varphi} ( \mathbb{H}^N) $ such that the corresponding solution to problem \eqref{Eqgeral-1} is global.
    \end{enumerate}
\end{theorem}

\begin{remark} The condition $\gamma<Q(p-1)$ in Theorem \ref{global-e-blowup-gamapositivo} (i) seems to be technical, but it appears naturally when choosing certain initial data, for example, if $u_0 \neq 0$ and $0\leq u_0(\eta) \leq C(1+|\eta|_{\mathbb{H}})^{-Q}$, then $u_0 \in L^{\infty}_{\varphi}( \mathbb{H}^N)$  if $\gamma<Q(p-1)$.
\end{remark}

The following result analyzes the case $\gamma < 0$.
\begin{theorem}\label{globalsolgamaneg} 
Assume $p>1$ and $-2< \gamma<0$. 
    \begin{enumerate}[(i)]
        \item If $1<p \leq  1 + {(2+\gamma)}/{Q},$
 then all nontrivial nonnegative solutions of problem \eqref{Eqgeral-1} blow up in finite time.

 \item If $ 1 + {(2+\gamma)}/{(Q+\gamma)} <p$, then there are nontrivial and nonnegative initial data in $L^\infty(\mathbb{H}^N)$  such that the corresponding solution to problem \eqref{Eqgeral-1} is global.
    \end{enumerate}

\end{theorem}

\begin{remark} Global solutions  are difficult to find when $\gamma < 0$ and $1 + (2+\gamma)/Q<p$. In \cite[Remark 3.2]{BIRINDELLI-A}, Birindelli and Capuzzo studied  the elliptic case of problem \eqref{Eqgeral-1} and found stationary solutions of the form $u(\eta) = C_{\varepsilon}(1 + |\eta|_{\mathbb{H}})^{-\frac{(Q - 2 - \varepsilon)}{2}}$, for $1 + {(2 + \gamma)}/{(Q - 2 - \varepsilon)} \leq p$. When $\mathbb{H}^N$ is replaced by $\mathbb{R}^N$, the Fujita critical exponent  for problem \eqref{Eqgeral-1} with $-2<\gamma < 0$ is $p_c=1 + (2 + \gamma)/N$ (see \cite{WANG}, \cite{QI}). We expected to find nontrivial global solutions for $1 + {(2 + \gamma)}/{Q} < p$ in item $(ii)$ of Theorem \ref{globalsolgamaneg}. However, this was not realizable due to some technical difficulties. One issue is the absence of comparison principles and classical solutions. Additionally, due to the lack of homogeneity in the Korányi norm, we were unable to obtain accurate estimates (without using Hölder's inequality) on the expression $e^{ t \Delta_{\mathbb{H}} } [ |\cdot |^{\gamma}_{\mathbb{H}} e^{\sigma \Delta_{\mathbb{H}} }u_0]$,
similar to the estimate obtained by Castillo et al. in the Euclidean frame, as noted in \cite[Lemma 4.1]{CASTILLO1}. 
This leaves the open problem of finding nontrivial global solutions for $$p \in \left( 1 + \frac{2 + \gamma}{Q}, 1 + \frac{2 + \gamma}{Q + \gamma} \right].$$
Similar problems exist in some specific domains of $\mathbb{R}^N$. For example, Bandle and Levine found difficulties in determining the Fujita critical exponent, as noted in \cite[Theorem 9.1 and Remark 2.4]{BandleLevine}, in the space 
$$
\mathbb{R}^N_k = \left\{ x = (x_1, \ldots, x_N) \in \mathbb{R}^N \;|\; x_i > 0, \; i = 1, \ldots, k \right\}
$$
with $ 0 < k \leq N$.  They established, considering $\mathbb{R}^N_k$ instead of $\mathbb{H}^N$, the following:
\begin{itemize}
\item If $1<p<1+ \frac{2+\gamma}{N+k}$ then any nonnegative and nontrivial solution of (\ref{Eqgeral-1}) can not be global.

\item $1+   \frac{2+\gamma}{N+k-2} <p< \left \{ \begin{array}{ll}
(N+1)/(N-3)& \mbox{ if }N\geq 4,\\
\mbox{ arbitrary} & \mbox{ if } N=2, 3\\
\end{array}\right.$ \\
~~\\
then problem (\ref{Eqgeral-1}) possesses global solutions.
\end{itemize}
Note that the result is not optimal due to the gap:
$$1+ \frac{2+\gamma}{N+k} <  1+ \frac{2+\gamma}{N+k-2}.$$ This issue was recently solved by the authors of the current article in \cite[2025]{CastilloHardy}.
\end{remark}

The remainder of the paper is organized as follows. In Section \ref{section2} we present the Heisenberg group with its properties and give some auxiliary results. The local existence of solutions to the problem (\ref{Eqgeral-1}) is proved in Section \ref{Exist}. Finally, Theorems \ref{global-e-blowup-gamapositivo} and \ref{globalsolgamaneg} are proved in Sections \ref{GloBl} and \ref{GloB2} respectively. 

\section{Preliminaries}\label{section2}

We present some concepts and results that will be useful in proving our findings.

\subsection{The Heisenberg group}
The Heisenberg group $\mathbb{H}^N$ is the space $\mathbb{R}^{2N+1} =\mathbb{R}^N \times \mathbb{R}^N \times \mathbb{R}$ equipped with the group operation
\begin{equation}\label{op.gr}
	(x,y, \tau ) \circ (\tilde{x}, \tilde{y}, \tilde{ \tau } ) = (x+ \tilde{x}, \,y+\tilde{y}, \,\tau +\tilde{ \tau } + 2  \left( x \cdot \tilde{y} - \tilde{x} \cdot y\right)),
\end{equation}
where $\cdot$ denotes the standard scalar product in $\mathbb{R}^N$. The identity element in $\mathbb{H}^N$ is $0$ and the inverse is $\eta^{-1} = -\eta$ for any $\eta \in \mathbb{H}^N$. There is a group of dilations on $\mathbb{H}^N$ given by
$$ \delta_r(x,y, \tau) = (rx,ry, r^2 \tau ) \quad \mbox{ for all } \,\, r>0.$$
A direct calculation shows that $\delta_r$ is an automorphism of $\mathbb{H}^N$ for every $r > 0$, and therefore $\mathbb{H}^N = (\mathbb{R}^{2N+1},\, \circ, \,\delta_r)$ is a homogeneous Lie group on $\mathbb{R}^{2N+1}$.

The left-invariant vector fields that span the Lie algebra $\mathfrak{h}_n$ of $\mathbb{H}^{N}$ are given by 
\begin{equation*}
	X_i = \partial_{x_i} + 2y_i \partial_{ \tau  } , \quad Y_i = \partial_{y_i} - 2x_i \partial_{ \tau  }, \quad \mathcal{T} = \partial_{\tau}, \quad i=1, \cdots, N .
\end{equation*}
Since $[X_i,Y_i] = - 4 \partial_{\tau}$ and $[X_i, X_j] = [Y_i, Y_j] = [X_i, Y_j] = 0 $ for all $i \neq j$, we have that the Lie algebra of $\mathbb{H}^N$ has stratification
$$ \mathfrak{h}_n = \mbox{span} \{X_1, \cdots, X_N, Y_1, \cdots, Y_N\} \oplus \mbox{span}\{\mathcal{T}\}.$$
Therefore, $\mathbb{H}^N$ is a $2$-step stratified Lie group with $2N$ generators. The Heisenberg gradient is $\nabla_{\mathbb{H}}=(X_1,\cdots,X_n, Y_1, \cdots, Y_n)$ and the sub-Laplacian on $\mathbb{H}^N$ (also referred to as Kohn Laplacian) is defined as
\begin{equation*}
	\Delta_{\mathbb{H}} = \sum_{i=1}^{N}  X_i^2 + Y_i^2 .
\end{equation*}
An explicit computation gives the  expression
\begin{equation*}
\begin{array}{ll}
\Delta_{\mathbb{H}}= \sum_{i=1}^{N} \left(  \partial^2_{x_ix_i} + \partial^2_{y_iy_i} \right) + 4(|x|^2+|y|^2) \partial_{\tau \tau }^2 +4 \sum_{i=1}^{N} \left ( y_i \partial^2_{x_i \tau} - x_i \partial^2_{y_i \tau}  \right). 
\end{array}
\end{equation*}
The operator $\Delta_{\mathbb{H}}$ is a degenerate elliptic operator. Since the Lie algebra of $\mathbb{H}^N$ fulfills the Hörmander condition, it is hypoelliptic (see \cite{Lars}). In sub-Riemannian geometry, stratified groups appear naturally. The Heisenberg group stands as a significant example of a sub-Riemannian manifold. Here, the concept of distance arises from a sub-Riemannian metric known as the Carnot-Carathéodory distance. This metric connects points in $\mathbb{H}^N$ through curves with tangent vectors in the span of $\{X_1, \ldots, X_N, Y_1, \ldots, Y_N\}$. The Carathéodory–Chow–Rashevsky Theorem ensures the existence of such curves connecting any two points in $\mathbb{H}^N$ (see \cite[Theorem 3.5.2]{Ovidiu}). A symmetric homogeneous norm in the group $ \mathbb{H}^N $ is a continuous function $ | \cdot |_{\mathbb{H}}: \mathbb{H}^N \rightarrow [0, +\infty)$ such that
\begin{itemize}
    \item $|\eta|_{\mathbb{H}} = 0$ if and only if $\eta=0,$
    \item $ |\delta_{\lambda} \eta|_{\mathbb{H}} = \lambda | \eta|_{\mathbb{H}}$ for all $\eta \in \mathbb{H}^N $, $\lambda>0,$
    \item $|\eta^{-1}|_{\mathbb{H}} = | \eta|_{\mathbb{H}} $ for all $\eta \in \mathbb{H}^N .$
\end{itemize}
It is established in \cite[Proposition 5.1.4]{BONFIGLIOLI} that any two homogeneous norms in $\mathbb{H}^N$ are equivalent. Moreover,  \cite[Theorem 5.2.8]{BONFIGLIOLI} shows that the Carnot-Carathéodory distance is a symmetric homogeneous norm on $\mathbb{H}^N$. Consequently, we may consider the following homogeneous norm on $\mathbb{H}^N$ (or an equivalent one):
\begin{equation}\label{koranyi-norm-def}
	|\eta|_{\mathbb{H}} = \left[  \left(  |x|^2 + |y|^2   \right)^2 + | \tau |^2  \right]^{\frac{1}{4}}.
\end{equation}
This gauge, $|\eta|_{\mathbb{H}}$, defined in \eqref{koranyi-norm-def} is commonly referred to as the Korányi norm, is equivalent to the Carnot-Carathéodory distance. Although it does not mirror the sub-Riemannian structure of $\mathbb{H}^N$, it proves to be very useful for calculations. Furthermore, $|\eta|_{\mathbb{H}}$ satisfies the triangle inequality $|\eta \circ \xi|_{\mathbb{H}} \leq |\eta|_{\mathbb{H}} + |\xi|_{\mathbb{H}}$, as seen in \cite[Proposition 5.1.8]{BONFIGLIOLI}. Through the group action, we can define a metric on $\mathbb{H}^N$ as follows: 
$$d(\eta, \xi) = |\xi^{-1} \circ \eta|_{\mathbb{H}}.$$
This metric induces a topology on $\mathbb{H}^N$. Thus, we can define the Heisenberg ball of $\mathbb{H}^N$, centered at $\eta$ and with radius $r > 0$, as  $B_{\mathbb{H}}(\eta, r) = \{ \xi \in \mathbb{H}^N \mid d(\eta, \xi) < r \}.$

In $\mathbb{R}^{2N+1}$, the Lebesgue measure remains invariant under the right-left translations of the Heisenberg group (see \cite[Proposition 1.3.21]{BONFIGLIOLI}). Consequently, since Haar measures on Lie groups are uniquely determined up to a constant factor, we represent the Haar measure on $\mathbb{H}^N$ by $d\varsigma$, which coincides with the  $2N+1$-dimensional Lebesgue measure. Additionally, the Lebesgue measure is homogeneous relative to the dilations $\delta_r$ for any $r>0$. This is expressed as, for a measurable set $E \subset \mathbb{R}^{2N+1}$ with Lebesgue measure $|E|$, we have
\begin{equation}\label{poly-growth}
    |\delta_r(E)| = r^Q |E|.
\end{equation}
Here, the exponent $Q := 2N+2$ denotes the homogeneous dimension of $\mathbb{H}^N$. The groups exhibiting the property described in \eqref{poly-growth} are known as groups of polynomial growth. Hence, the Heisenberg group is identified as a 2-step stratified Lie group with polynomial growth, where the homogeneous dimension is given by $Q = 2N + 2$.

The operator $\Delta_{\mathbb{H}}$ is a self-adjoint operator on $L^2(\mathbb{H}^N)$ (see \cite[Proposition 2.1]{HULANICKI1}). By \cite[Theorem 3.4]{HULANICKI1} exists a unique symmetric Markov semigroup of
operators $\{e^{t\Delta_{\mathbb{H}}}\}_{t\geq 0}$ on $L^p(\mathbb{H}^N)$, $1 \leq p \leq \infty$, associated with $\Delta_{\mathbb{H}}$, such that for each $p$ and each $t > 0$, $  e^{t\Delta_{\mathbb{H}}} $ is a bounded operator on $L^p(\mathbb{H}^N)$. It is well known that the semigroup $e^{t\Delta_{\mathbb{H}}}$ admits a kernel $h_t$, explicitly given by
$$  h_t (x, \tau) = \dfrac{1}{(2\pi)^{N+2} 2^N} \int_{\mathbb{R}} \left( \dfrac{\lambda}{\sinh t \lambda }\right)^N \exp \left\{ \dfrac{-|x|^2 \lambda}{ 4 \tanh t \lambda } \right\} e^{i \lambda \tau } d\lambda , \quad (x,\tau) \in \mathbb{R}^{2N }\times \mathbb{R}$$
(see \cite[Theorem 1.3.2]{RANDALL}), and the action of the semigroup is given by
$$ e^{t\Delta_{\mathbb{H}}} u (\eta) = ( u \ast h_t ) (\eta)= \int_{\mathbb{H}^N} u (\varsigma) h_t (\varsigma^{-1} \circ \eta) d\varsigma =   \int_{\mathbb{H}^N} u ( \eta \circ \varsigma^{-1}) h_t (\varsigma) d\varsigma. $$
The kernel $h_t(\eta)$ is a $C^{\infty}$ function in the variable $(t,\eta) \in (0, \infty) \times \mathbb{H}^N$, $\|h_t\|_{L^1(\mathbb{H}^N)} = 1$  for any $t>0$ and satisfies the heat equation $\partial_t - \Delta_{\mathbb{H}} = 0$ (see \cite{VAROPOULOS}, \cite{HULANICKI1} ). Furthermore, there exist two positive constants $c_\star, C_\star$ such that the kernel can be estimated as follows:
\begin{equation}\label{estimativadokernel}
	c_{\star}t^{-\frac{Q}{2}} \exp \left(  - \dfrac{C_\star|\eta|_{\mathbb{H}}^2}{t} \right) \leq h_t(\eta) \leq C_\star t^{-\frac{Q}{2}} \exp \left(  - \dfrac{c_\star|\eta|_{\mathbb{H}}^2}{t} \right) 
\end{equation}
for any $t > 0$ and $\eta \in \mathbb{H}^N$ (see \cite[Theorems IV. 4.2 and IV.4.3]{VAROPOULOS}).

The semigroup $\{e^{t\Delta_{\mathbb{H}}}\}_{t>0}$ is a contraction semigroup on $L^p(\mathbb{H}^N)$, $1 \leq p\leq \infty$, which is strongly continuous for $p< \infty $. Due to Stein's interpolation theorem, one can show that $\{e^{t\Delta_{\mathbb{H}}}\}_{t>0}$ is an analytic semigroup in $L^p(\mathbb{H}^N)$, $p< \infty $ (see \cite[Theorem 1.4.2 ]{DAVIES1} or \cite[Theorem 3.1]{FOLLAND1}). Furthermore, since $\mathbb{H}^N$ is a Lie group with polynomial growth, the semigroup $\{e^{t\Delta_{\mathbb{H}}}\}_{t>0}$ is bounded and analytic in $L^1( \mathbb{H}^N)$ (see \cite[Corollary VIII. 2.6]{VAROPOULOS}). 

\subsection{Auxiliary results} 
\begin{lemma}[Reverse H\"older inequality]\label{invholder} Suppose $p \in (1, \infty)$, $f\in L^{1/p}(\mathbb{H}^N)$  and 
	$$0<\int_{\mathbb{H}^N} |g(\varsigma)|^{-\frac{1}{p-1}}d\varsigma<\infty.$$
	Then 
	$$	\int_{\mathbb{H}^N}  |f(\varsigma)g(\varsigma)| d\varsigma \geq \left( \int_{\mathbb{H}^N} |f(\varsigma)|^{\frac{1}{p}}d\varsigma \right)^p \left( \int_{\mathbb{H}^N} |g(\varsigma)|^{-\frac{1}{p-1}} d\varsigma \right)^{-(p-1)}.$$
\end{lemma}
For a proof, see \cite[Theorem 2.12]{Adams}. 

Let $E$ be a measurable set of finite measures in $\mathbb{R}^N$, i.e., $|E|<\infty$. The symmetric rearrangement of $E$, which we denote by $E^*$, is the open ball centered at the origin whose measure coincides with the measure of $E$. Let $f : \mathbb{R}^N \to [0, \infty)$ be a measurable function vanishing at infinity, that is, $|\{x; |f(x)|>t\}|$ is finite for all $t>0$. The symmetric decreasing rearrangement of $f$, denoted by $f^*$, is defined as  
$$
{f^*}(x) = \int_0^{\infty} \chi_{\{y;\,|f(y)|>t\}^*} (x) \, dt.
$$
See more details in \cite[Chapter 3]{Lieb1} and \cite[Chapter X]{inequalities}.

We now state the Riesz's rearrangement inequality, a classical tool in functional analysis and a key result in the study of sharp inequalities, which will be essential for proof of the Lemma \ref{Lemachave1}.

\begin{lemma}[Riesz’s rearrangement inequality] Let $f, g, h: \mathbb{R}^N \rightarrow \mathbb{R}^+$ be measurable functions that vanish at infinity. Then they satisfy the inequality:
\begin{equation} \label{Riesz1}
\int_{\mathbb{R}^N} \int_{\mathbb{R}^N} f(z) g(z-y) h(y) \, dz\,dy
\leq \int_{\mathbb{R}^N} \int_{\mathbb{R}^N} f^{\star}(z) g^{\star}(z-y) h^{\star}(y) \, dz\,dy,
\end{equation}
in the sense that the left hand side is finite whenever the right hand side is finite, where $f^\star : \mathbb{R}^N \rightarrow \mathbb{R}^+$, $g^\star : \mathbb{R}^N \rightarrow \mathbb{R}^+ $ and $h^\star : \mathbb{R}^N \rightarrow \mathbb{R}^+$ are the symmetric decreasing rearrangement of the functions $f$, $g$ and $h$ respectively.
\end{lemma}
{\bf Proof.}
    See \cite[Theorem 3.7]{Lieb1}
 \cqd \\

The next result has already been proved  in \cite[Lemma 6]{PINSKY1} for $k=2$. We extend it here for all $k>0$. 
\begin{lemma}\label{Lemachave1} For $m \leq 0$, $c, k>0$, $\delta \geq 0$ ($-N<m\leq 0$ if $\delta=0$) and $t>0$  the function $G: \mathbb{R}^N \to \mathbb{R}$ defined by 
$$G(x) = \int_{\mathbb{R}^N} \exp\left( -c \dfrac{|x-y|^{k}}{t} \right) (\delta+|y|)^m dy$$
attains its maximum at $x=0$.
\end{lemma}
{\bf Proof.} Fix $x \in \mathbb{R}^N$ and for each $\varepsilon \in (0,1)$ define $\psi_{x, \varepsilon}= |B(0,\varepsilon)|^{-1} \; \chi_{B(x,\varepsilon)}$. Note that the symmetric decreasing rearrangement is $\psi_{x,\varepsilon}^\star=|B(0, \varepsilon)|^{-1}\;  \chi_{B(0,\varepsilon)}=\psi_{0,\varepsilon}$.

Using Riesz's rearrangement inequality (\ref{Riesz1}) by setting $g(z-y) = \exp (-c|z-y|^k/t)$, $h(y)=(\delta+|y|)^m$ and $f(z)= \psi_{\varepsilon}(z)$, we obtain
\begin{eqnarray*}
 &\int_{\mathbb{R}^N} \int_{\mathbb{R}^N} \psi_{x,\varepsilon}(z) g(z-y) h(y) dz dy \\
 \noalign{\medskip}
 &\leq \int_{\mathbb{R}^N} \int_{\mathbb{R}^N}  \psi_{0,\varepsilon}(z) g(z-y) h(y) dz dy.
\end{eqnarray*}
Observe that the integral on the right-hand side is finite; therefore, by Fubini's theorem, we obtain
\begin{equation}\label{Rea.un}
\begin{array}{ll}
 &\int_{\mathbb{R}^N} \left [ \int_{\mathbb{R}^N} \psi_{x,\varepsilon}(z) g(z-y)  dz \right] h(y) dy \\
 \noalign{\medskip}
 &\leq \int_{\mathbb{R}^N} \left [ \int_{\mathbb{R}^N}  \psi_{0,\varepsilon}(z) g(z-y) dz \right] h(y) dy.
\end{array}
\end{equation}

Define $\Phi_{x,\varepsilon}: \mathbb{R}^N \to \mathbb{R}$ by 
$$
\begin{array}{rll}
\Phi_{x,\varepsilon}(y) &= \int_{\mathbb{R}^N} \psi_{x,\varepsilon}(z) g(z-y) dz.
\end{array}
$$
Then, $\Phi_{x,\varepsilon}(y) \to g(x-y)$ a.e. $y \in \mathbb{R}^N$ as $\varepsilon \to 0$. Moreover, for $|y|\geq 2|x|+2$ and $|z-x|\leq \varepsilon$ we have
$$|z-y|\geq |y-x|-|z-x|\geq |y-x|-\epsilon\geq |y|-|x|-\epsilon \geq \frac{|y|}{2},$$
and
\begin{equation}\label{Rea.do}
\begin{array}{ll}
|\Phi_{x,\varepsilon}(y)|&\leq |B(0,\varepsilon)|^{-1}\int_{|z-x|\leq \varepsilon} \exp \left (-c\frac{|z-y|^k}{t}\right ) dz \\
& \leq \exp \left(-c\frac{|y|^k}{2^k t} \right).
\end{array}
\end{equation}
Let $\Psi: \mathbb{R}^N \to \mathbb{R}$ given by 
$$
\Psi(z)=\left \{ \begin{array}{ll}
h(z)& \mbox{ if } |z|\leq 2|x|+2, \\
\exp \left(-c\frac{|z|^k}{2^k t} \right) h(z)& \mbox{ if } |z|\geq 2|x|+2.
\end{array}\right.
$$
From (\ref{Rea.do}) we see that $|\Phi_{x,\varepsilon}(y) h(y)|\leq \Psi(y)$ for all $\varepsilon \in (0, 1)$. Since $\Psi \in L^1(\mathbb{R}^N)$, by Lebesgue's dominated convergence theorem, letting $\varepsilon \to 0$ in the inequality (\ref{Rea.un}), we obtain the desired result.

 The limit case $\delta = 0$ is handled as in \cite[Lemma 2.10]{CASTILLO1}. In fact, this follows from the case $\delta>0$ and the monotone convergence theorem, letting $\delta \to 0^+$, since $-N<m \leq 0$. 
 \cqd \\

\begin{remark}\label{eq-das-normas} According to \cite[Proposition 5.1.4]{BONFIGLIOLI}, any two homogeneous norms in $\mathbb{H}^N$ are equivalent. This means that the Korányi norm defined in \eqref{koranyi-norm-def} can be replaced with a more straightforward homogeneous norm when necessary to facilitate the calculations. This was used, for instance, in \cite[Proposition 4.1]{PALMIERI}.
\end{remark}

\begin{lemma}\label{UpperSemigroup} Let $u_0 \in L^\infty(\mathbb{H}^{N})$ and  $0 < \gamma < Q$. Then there exists a constant $C>0$ depending on $ \| u_0\|_{L^\infty(\mathbb{H}^{N})},$ $N$ and $\gamma$ such that
\begin{equation} \label{Upper1} 
\left \| t^{-\frac{Q}{2}}\int_{\mathbb{H}^N}  \exp\left(  -c_\star \dfrac{| \varsigma^{-1} \circ \eta |^{2}_{\mathbb{H} }}{t}\right) \, | \varsigma |_{\mathbb{H}}^{-\gamma} \, u_0(\varsigma) \, d\varsigma  \right \|_{L^\infty(\mathbb{H}^{N})} \leq C \, t^{-\frac{\gamma}{2}},
\end{equation}
for all $t>0.$
\end{lemma}
{\bf Proof.} Note that by \cite[Corollary 5.4.5]{BONFIGLIOLI}, the function $\varsigma \mapsto |\varsigma|_{\mathbb{H}}^{-\gamma}$ is locally integrable in $\mathbb{H}^N$ if $0<\gamma< Q$. 

Recall that the natural Haar measure in $\mathbb{H}^N $ is identical to the Lebesgue measure $d\varsigma = dx' dy'd \tau'$ in $\mathbb{R}^{2N+1} =\mathbb{R}^N \times \mathbb{R}^N \times \mathbb{R}$. As it was shown in \cite[Proposition 4.1]{PALMIERI}, we can argue that
\begin{align} \label{integral do ptofixo}
\begin{aligned}
	&t^{-\frac{Q}{2}}\int_{\mathbb{H}^N}  \exp\left(  -c_\star \dfrac{| \varsigma^{-1} \circ \eta |^{2}_{\mathbb{H} }}{t}\right) \, | \varsigma |_{\mathbb{H}}^{-\gamma} \, u_0(\varsigma) \, d\varsigma  \\
\noalign{\medskip}
	 &\leq C \, t^{-\frac{Q}{2}}	\int_{\mathbb{R}^{2N+1}}  \exp\left(  -c_1 \dfrac{  |x-x' |^2 + |y - y'|^2 + | \tau - \tau' + 2 (x'\cdot y - x \cdot y' ) | }{t}\right) \times \\
	& \ \ \ \ \ \ \ \  \  \  \ \left( |x'| +|y'|+|\tau'|^{\frac{1}{2}} \right)^{-\gamma } dx' dy'd \tau' ,
\end{aligned}
\end{align}
where $\eta = (x,y, \tau)$, $\varsigma = (x',y', \tau ')$ and $C, c_1$ are positive constants. Here, we have used Remark \ref{eq-das-normas} and replaced the Korányi norm with a more straightforward homogeneous norm. Thus, using Young's inequality, Fubini's theorem and Lemma \ref{Lemachave1}, we have
\begin{align}\label{use-cal-fren}
    \begin{aligned}
      &  t^{-\frac{Q}{2}} \int_{\mathbb{H}^N}  \exp\left(  -c_\star \dfrac{| \varsigma^{-1} \circ \eta |^2_{\mathbb{H}}}{t}\right)| \varsigma |_{\mathbb{H}}^{-\gamma}  \, u_0(\varsigma) \, d\varsigma \\
      \noalign{\medskip}
	&  \leq  C \, t^{-\frac{Q}{2}}	\int_{\mathbb{R}^{2N+1}}  \exp\left(  -c_1 \dfrac{  |x-x' |^2 + |y - y'|^2 }{t}\right) 
	\exp \left(  \dfrac{-c_1 | \tau - \tau' + 2 (x'\cdot y - x \cdot y' ) | }{t}  \right) \times\\
 & \ \ \ \ \ \ \ \  \  \  \left(|x'| +|y'|+|\tau'|^{\frac{1}{2}} \right)^{-\gamma } dx' dy'd \tau' \\
 \noalign{\medskip}
	&  \leq  C \, t^{-\frac{Q}{2}}	\int_{\mathbb{R}^{2N+1}}  \exp\left(  -c_1 \dfrac{  |x-x' |^2 + |y - y'|^2 }{t}\right)  \exp \left(  \dfrac{-c_1 | \tau - \tau' + 2 (x'\cdot y - x \cdot y' ) | }{t}  \right) \times\\
 & \ \ \ \ \ \ \ \  \  \ \left( |x'|^{\frac{N}{Q}}  |y'|^{\frac{N}{Q}}  |\tau'|^{\frac{1}{Q}}    \right)^{-\gamma }    dx' dy'd \tau' \\
 \noalign{\medskip}
	& \leq  C \, t^{-\frac{Q}{2}} \int_{\mathbb{R}^N}   \exp\left(  -c_1 \dfrac{  |x-x' |^2 }{t}\right)  |x '|^{\frac{-\gamma N}{Q}} \int_{\mathbb{R}^N}  \exp\left(  -c_1 \dfrac{ |y - y'|^2 }{t}\right) |y '|^{-\frac{\gamma N}{Q}} \times \\
	& \ \ \ \ \ \ \ \  \  \  \underbrace{\int_{\mathbb{R}}   \exp \left(  \dfrac{-c_1 | \tau - [\tau' - 2 (x'\cdot y - x \cdot y' )] | }{t}  \right)  |\tau '|^{-\frac{\gamma}{Q}}     d \tau'}_{{I}} \, dy' dx'. \\
    \end{aligned}
\end{align}
Now, we will estimate the integral $I$. From Lemma \ref{Lemachave1}, we obtain
\begin{align*}
    \begin{aligned}
 I &= \int_{\mathbb{R}}   \exp \left(  \dfrac{-c_1 \left | \, [\tau +  2 (x'\cdot y - x \cdot y' )] - \tau' \, \right| }{t}  \right)  |\tau '|^{-\frac{\gamma}{Q}}     d \tau' \\ 
\noalign{\medskip}
& \leq \int_{\mathbb{R}} \exp \left(  \dfrac{-c_1 |\tau' | }{t}  \right)  |\tau '|^{-\frac{\gamma}{Q}}   d \tau'.
\end{aligned}
\end{align*}
Thus, from this and  \eqref{use-cal-fren}, we have 
\begin{eqnarray*}
 & & t^{-\frac{Q}{2}} \int_{\mathbb{H}^N}  \exp\left(  -c_\star \dfrac{| \varsigma^{-1} \circ \eta |^2_{\mathbb{H}}}{t}\right)| \varsigma |_{\mathbb{H}}^{-\gamma} \, u_0(\varsigma) \,  d\varsigma \\
 &\leq &  C \, t^{-\frac{Q}{2}} \int_{\mathbb{R}^N}   \exp\left(  -c_1 \dfrac{  |x' |^2 }{t}\right)  |x '|^{\frac{-\gamma N}{Q}} d x' \\
 & & \times \int_{\mathbb{R}^N}  \exp\left(  -c_1 \dfrac{ |y'|^2 }{t}\right) |y '|^{-\frac{\gamma N}{Q}} d y' \\
 & & \times \int_{\mathbb{R}} \exp \left(  \dfrac{-c_1 |\tau' | }{t}  \right)  |\tau '|^{-\frac{\gamma}{Q}}   d \tau'.
\end{eqnarray*}
Hence, applying the changes of variables $x' = t^{\frac{1}{2}}x$, $y' = t^{\frac{1}{2}}y$, and $\tau' = t \tau$, we can arrive at the following conclusion
\begin{equation*} 
t^{-\frac{Q}{2}} \int_{\mathbb{H}^N}  \exp\left(  -c_\star \dfrac{| \varsigma^{-1} \circ \eta |^2_{\mathbb{H}}}{t}\right)| \varsigma |_{\mathbb{H}}^{-\gamma}  \, u_0(\varsigma) \, d\varsigma \leq C \, t^{-\frac{Q}{2}} t^{\frac{Q}{2}} t^{-\gamma \left(  \frac{N}{2Q}+\frac{N}{2Q} +\frac{2}{2Q} \right)}  = C \, t^{-\frac{\gamma}{2}}.      
\end{equation*}
The proof is now complete.

\cqd

\section{Local existence} \label{Exist}

The local existence for $\gamma\geq 0$ is given in the following result.
\begin{theorem}\label{existenciagamapos} Assume $p>1, \gamma\geq 0$,  $u_0 \in L^{\infty}_{\varphi}(\mathbb{H}^N),$ $u_0\geq 0$, $u_0 \neq 0$ and $\varphi(\eta) = (1+| \eta|_{\mathbb{H}})^{\frac{\gamma}{p-1}}$. There exist  $T > 0$  and a unique function $u \in L^\infty ([0,T], L^{\infty}_{\varphi}(\mathbb{H}^N))$ solution of problem \eqref{Eqgeral-1}. Moreover, it is possible to extend this solution to a maximal interval $(0,T_{\max})$ so that we have the blow up alternative:  either $T_{\max}= +\infty$ or  $T_{\max}<\infty$ and  $\lim_{t \rightarrow T_{\max}} \| \varphi \cdot u(t)\|_{L^{\infty}(\mathbb{H}^N)} = + \infty$.
\end{theorem}

{\bf Proof. } According to \cite[Theorem 3.1]{FOLLAND1} we have
$h_t (\xi) = \lambda^Q h_{\lambda^2 t} (\delta_{\lambda}  (\xi) ), \quad t>0,$
for any $\lambda>0$ and any $\xi \in \mathbb{H}^N$, where $\delta_{\lambda}$ is the anisotropic dilation (defined in Section \ref{section2}). In particular, when $\lambda = t^{ - \frac{1}{2}}$ we have
\begin{equation}\label{invarianciaheat}
	h_t (\xi) = t^{- \frac{Q}{2}} h_{1} (\delta_{  t^{ - {1}/{2}}  }  (\xi) )
\end{equation}
for any $\xi \in  \mathbb{H}^N$ and any $t>0.$

We use the Banach fixed point theorem for the mapping $\mathcal{F}: B_K \rightarrow X$ defined by 
\begin{equation}\label{mild-fixed-point-use}
	\mathcal{F}u(t)= e^{ t\Delta_{\mathbb{H}} }u_0 + \underbrace{  \int_{0}^{t} e^{(t-s)\Delta_{\mathbb{H}}}| \cdot |_{\mathbb{H}}^{\gamma} u^p(s) ds }_{\mathcal{G}u(t)},
\end{equation}
where  
$$X = L^\infty ([0,T], L^{\infty}_{\varphi}(\mathbb{H}^N)) =  \left\{ v: [0,T] \rightarrow  L^{\infty}_{\varphi}(\mathbb{H}^N) ; \sup_{t \in [0,T] }\| \varphi  \cdot v (t)\|_{L^{\infty} ( \mathbb{H}^N )} < \infty \right \} ,$$ 
and 
$B_K = B_K^T$ denotes the closed ball in $X$ with center at $0$ and radius $K$. The positive constants $K$ and $T$ will be chosen so that $\mathcal{F}$ satisfies $\mathcal{F}: B_K \rightarrow B_K$ and is a strict contraction. We start by estimating $ \| \varphi \cdot e^{ t\Delta_{\mathbb{H}} }u_0 \|_{L^{\infty} ( \mathbb{H}^N )} $. Note that
\begin{align}\label{est.1pf-pos}
	\begin{aligned}
		[e^{ t\Delta_{\mathbb{H}} }u_0](\eta) = \int_{\mathbb{H}^N} h_t (\xi) u_0( \eta \circ \xi^{-1} ) d\xi = t^{- \frac{Q}{2}} \int_{\mathbb{H}^N} h_1 (\delta_{  t^{ - {1}/{2}}  }  (\xi) ) u_0( \eta \circ \xi^{-1} ) d\xi.
	\end{aligned}
\end{align}
By applying the change of variables $\zeta = \delta_{t^{-{1}/{2}}}(\xi)$ to (\ref{est.1pf-pos}),  which is an isomorphism on $\mathbb{H}^N$ with $\delta_r^{-1} = \delta_{r^{-1}}$, we obtain
$$  [e^{ t\Delta_{\mathbb{H}} }u_0](\eta) = \int_{\mathbb{H}^N} h_1 (\zeta)  u_0( \eta \circ (   \delta_{  t^{ {1}/{2}}  }  (\zeta))^{-1} ) d \zeta.$$
By the estimate (\ref{estimativadokernel}) 
\begin{align}\label{est.2pf-pos}
\begin{aligned}
		[e^{ t\Delta_{\mathbb{H}} }u_0](\eta) & \leq  \| \varphi \cdot u_0 \|_{L^{\infty} ( \mathbb{H}^N )}  \int_{\mathbb{H}^N} h_1 (\zeta)(1+| \eta \circ (   \delta_{  t^{ {1}/{2}}  }  (\zeta))^{-1} |_{\mathbb{H}})^{\frac{-\gamma}{p-1}}   d \zeta\\
  \noalign{\medskip}
		&  \leq \| \varphi \cdot u_0 \|_{L^{\infty} ( \mathbb{H}^N )}  \left(  \int_{R_1} + \int_{R_2} \right) h_1 (\zeta)(1+| \eta \circ (   \delta_{  t^{ {1}/{2}}  }  (\zeta))^{-1} |_{\mathbb{H}})^{\frac{-\gamma}{p-1}} d \zeta\\
		&  \mbox{ where } R_1 = \left \{ \zeta  ; \,   | \eta \circ (   \delta_{  t^{ {1}/{2}}  }  (\zeta))^{-1} |_{\mathbb{H}} \leq | \eta |_{\mathbb{H}}/2 \, \right \}  \mbox{ and } \\
            &  R_2 = \left\{ \zeta  ; \,   | \eta \circ (   \delta_{  t^{ {1}/{2}}  }  (\zeta))^{-1} |_{\mathbb{H}} \geq | \eta |_{\mathbb{H}}/2 \, \right \}\\
  \noalign{\medskip}
	&  \leq  C\| \varphi \cdot  u_0 \|_{L^{\infty} ( \mathbb{H}^N )} \left[  \int_{R_1} e^{-c_\star |\zeta |_{\mathbb{H}}^2} d\zeta + \left(  1 + \dfrac{|\eta|_{\mathbb{H}}}{2}\right)^{\frac{-\gamma}{p-1}} \int_{\mathbb{H}^N}   e^{-c_\star |\zeta |_{\mathbb{H}}^2} d\zeta       \right].
\end{aligned}
\end{align}
Using triangle inequality on the region $R_1$, we have
\begin{equation} \label{Triang1}
\dfrac{| \eta |_{\mathbb{H}}}{2} \geq | \eta \circ ( \delta_{ t^{ {1}/{2}} } (\zeta))^{-1} |_{\mathbb{H}} \geq |\eta|_{\mathbb{ H}} - t^{\frac{1}{2}}| \zeta |_{\mathbb{H}}.
\end{equation}
So, $R_1 \subset \left\{ \zeta ; | \zeta |_{\mathbb{H}} \geq t^{ - \frac{1}{2}} {| \eta |_{\mathbb{H}}}/{2} \right\}$. Using \eqref{Triang1} and the fact that the Haar measure is similar to the Lebesgue measure, we can continue computing \eqref{est.2pf-pos} as follows:
\begin{align}\label{est.3pf-pos}
	\begin{aligned}
		[e^{ t\Delta_{\mathbb{H}} }u_0](\eta) 
		&  \leq C\| \varphi \cdot u_0 \|_{L^{\infty} ( \mathbb{H}^N )} \left[  \int_{R_1} e^{-c_\star |\zeta |_{\mathbb{H}}^2} d\zeta + \left(  1 + \dfrac{|\eta|_{\mathbb{H}}}{2}\right)^{\frac{-\gamma}{p-1}} \int_{\mathbb{H}^N}   e^{-c_\star |\zeta |_{\mathbb{H}}^2} d\zeta       \right]\\
  \noalign{\medskip}
  &  \leq C\| \varphi \cdot u_0 \|_{L^{\infty} ( \mathbb{H}^N )} \left[ e^{\frac{-c_\star| \eta |^2_{\mathbb{H}}}{8t}} \int_{ \mathbb{H}^N  } e^{-c_\star \frac{ |\zeta |_{\mathbb{H}}^2}{2}} d\zeta + C_1\left(  1 + \dfrac{|\eta|_{\mathbb{H}}}{2}\right)^{\frac{-\gamma}{p-1}}     \right]\\
  \noalign{\medskip}
&	\leq C\| \varphi \cdot u_0 \|_{L^{\infty} ( \mathbb{H}^N )} \left(  e^{\frac{-c_\star| \eta |^2_{\mathbb{H}}}{8t}} \int_{\mathbb{R}^N} \int_{\mathbb{R}^N}  \int_{ -\infty }^{+\infty} e^{-c_1 \frac{|x |^2}{2}} e^{-c_1 \frac{|y |^2}{2}} e^{-c_1 \frac{|\tau|}{2} } d\tau dy dx + C_1 (\varphi(\eta))^{-1}      \right)\\
\noalign{\medskip}
& \leq C\| \varphi \cdot u_0 \|_{L^{\infty} ( \mathbb{H}^N )} \left(   \varphi(\eta) e^{\frac{-c_\star| \eta |^2_{\mathbb{H}}}{8}} C_2  + C_1   \right) (\varphi(\eta))^{-1} \\
\noalign{\medskip}
& \leq C \| \varphi \cdot u_0 \|_{L^{\infty} ( \mathbb{H}^N )} (\varphi(\eta))^{-1},
	\end{aligned}
\end{align}
for $0 < t<T < 1$ (note that the function $\eta \mapsto \varphi(\eta)e^{\frac{-c_\star| \eta |^2_{\mathbb{H}}}{8}}$ is bounded). Thus, we deduce that
\begin{equation}\label{use-more-n=one}
   \| \varphi \cdot  e^{ t\Delta_{\mathbb{H}} }u_0 \|_{L^{\infty} ( \mathbb{H}^N )} \leq C \| \varphi \cdot u_0 \|_{L^{\infty} ( \mathbb{H}^N )},  
\end{equation}
for $t \in [0,T]$, $T<1.$  

Now, we will estimate $\mathcal{G}u(t)$ arguing similarly as in the derivation of estimates $\eqref{est.2pf-pos}$ and $\eqref{est.3pf-pos}$. Indeed,
\begin{align}\label{est.4pf-pos}
	\begin{aligned}
		[\mathcal{G}u(t)](\eta) & =\int_{0}^{t} \left [e^{(t-s)\Delta_{\mathbb{H}}}| \cdot |_{\mathbb{H}}^{\gamma} u^p(s) \right ](\eta) \, ds\\
  \noalign{\medskip}
		&  =	\int_{0}^{t}  \int_{\mathbb{H}^N} h_1 (\zeta) \, \left | \eta \circ (   \delta_{  (t-s)^{ {1}/{2}}  }  (\zeta))^{-1} \right |_{\mathbb{H}}^{\gamma} \, u^p( \eta \circ (   \delta_{  (t-s)^{ {1}/{2}}  }  (\zeta))^{-1} , s) \, \, d \zeta ds\\
  \noalign{\medskip}
		&	\leq C\int_{0}^{t}  \| \varphi \cdot u (s) \|^p_{L^{\infty} ( \mathbb{H}^N )}   \int_{\mathbb{H}^N} h_1 (\zeta) \, \, \left | \eta \circ (   \delta_{  (t-s)^{ {1}/{2}}  }  (\zeta))^{-1} \right |_{\mathbb{H}}^{\gamma} \\
  &  \,\,\,\,\,\,\,\,\,\,\,\,\,\, \left ( 1+ \left | \eta \circ (   \delta_{  (t-s)^{ {1}/{2}}  }  (\zeta))^{-1} \right |_{\mathbb{H}} \right )^{\frac{-\gamma p}{p-1}}  d \zeta   ds   \\
  \noalign{\medskip}
		& \leq C K^p\int_{0}^{t} \int_{\mathbb{H}^N} h_1 (\zeta) \, \left( 1+ \left | \eta \circ (   \delta_{  (t-s)^{ {1}/{2}}  }  (\zeta))^{-1} \right |_{\mathbb{H}} \right)^{\frac{-\gamma }{p-1}}  d \zeta   ds   \\
  \noalign{\medskip}
		& \leq C K^p\int_0^t \left(   e^{-\frac{c_\star| \eta |^2_{\mathbb{H}}}{8}} C_2 + C_1 \varphi^{-1}(\eta)   \right)ds\\
  \noalign{\medskip}
		& \leq C K^p T \varphi^{-1}(\eta).
	\end{aligned}
\end{align}
Based on the estimates \eqref{use-more-n=one} and \eqref{est.4pf-pos}, there is a constant $C_3 >0 $ such that
\begin{equation}\label{Hen1}
\left\| \varphi \cdot \mathcal{F} u(t) \right \|_ {L^{\infty}(\mathbb{H}^N)} \leq C_3 \left ( \| \varphi \cdot u_0 \|_{L^{\infty} ( \mathbb{H}^N )} + K^p \, T \right ) .
\end{equation}
Similarly, for $u, v \in B_K$, we have
\begin{align}\label{Hen2}
	\begin{aligned}
		\| \varphi \cdot	(\mathcal{F}u (t) - 	\mathcal{F}v (t)) \|_{L^{\infty}(\mathbb{H}^N)} \leq  C_4 K^{p-1}  T 	\sup_{t \in [0,T] } \| 	\varphi \cdot  (u (t) - 	v (t)) \|_ {L^{\infty}(\mathbb{H}^N)}
	\end{aligned}
\end{align}
for some constant $ C_4 >0$. Thus, choosing 
$K>  2 C_3 \| \varphi \cdot u_0 \|_{L^{\infty}(\mathbb{H}^N)}$ 
and $0<T<1$ sufficiently small so that
$\max \{C_3, \,C_4\} K^{p-1} T < 1/2,$
 from \eqref{Hen1} and \eqref{Hen2}, we can conclude that $\mathcal{F}$ is a contraction in $B_K$, and it has a unique fixed point $u$ in $B_K$, which is a local mild solution of \eqref{Eqgeral-1}.

The uniqueness in $X$ follows from the following: if $u,v$ are solutions on $[0,T]$ of \eqref{mild-fixed-point-use}, then by performing a calculation similar to \eqref{est.4pf-pos}, we obtain
\begin{align*}
   \| \varphi \cdot	(u (t) - 	v (t)) \|_{L^{\infty}(\mathbb{H}^N)} &\leq C \left(\sup_{s \in [0,T]} \| \varphi \cdot u(s) \|_{L^{\infty}(\mathbb{H}^N)}^{p-1} +  \sup_{s \in [0,T]} \| \varphi \cdot v(s) \|_{L^{\infty}(\mathbb{H}^N)}^{p-1} \right) \times  \\
   & \qquad \qquad \int_{0}^{t} \| 	\varphi \cdot (u (s) - 	v (s)) \|_{L^{\infty}(\mathbb{H}^N)} d s. 
\end{align*}
Applying Gronwall’s inequality, it follows that $u=v$. 

Since $T$ depends on $\|\varphi \cdot u_0\|_ {L^{\infty}(\mathbb{H}^N)}$ the blow-up alternative follows by standard arguments (see \cite[Remark 2.4.2 and Theorem 2.4.3]{breziscaz}).

\cqd

\medskip \noindent
The local existence for $\gamma<0$ is given in the following result.
\begin{theorem}\label{existenciagamneg} Suppose that $p>1$, $N\geq 1$ , $-2< \gamma<0$, $u_0\in L^\infty( \mathbb{H}^N ), u_0\geq 0 $ and $u_0 \neq 0$. There exist $T > 0$ and a unique function $u \in L^{\infty}([0,T] ; L^\infty(\mathbb{H}^N))$ solution to the problem \eqref{Eqgeral-1}. Moreover, it is possible to extend this solution to a maximal interval $(0,T_{\max})$ so that we have the blow-up alternative:  either $T_{\max}= +\infty$ or  $T_{\max}<\infty$ and  $\lim_{t \rightarrow T_{\max}} \| u(t)\|_{L^{\infty}(\mathbb{H}^N)} = + \infty$.
\end{theorem}
{\bf Proof. }  Let us determine the solution in the Banach space 
$$X = L^{\infty}([0,T] ; L^\infty(\mathbb{H}^N)) = \left\{\phi:[0,T] \rightarrow L^\infty(\mathbb{H}^N) ; \sup_{t \in [0,T]} \| \phi (t)\|_{L^{\infty}(\mathbb{H}^N)} < \infty \right \}.$$
Let $B_K = B_K^T$ denote the closed ball in $X$ with center $0$ and radius $K$, where $K>0$ and $T>0$ will be chosen so that $\mathcal{F}$ defined in \eqref{mild-fixed-point-use} satisfies $\mathcal{F}: B_K \rightarrow B_K$ and is a strict contraction. Using the estimate \eqref{estimativadokernel} we have 
$$
\begin{array}{ll}
	[\mathcal{G}u(t)](\eta)  &\leq   C\int_{0}^{t} \int_{\mathbb{H}^N} (t-s)^{-\frac{Q}{2}} \exp\left(  -c_\star \dfrac{| \varsigma^{-1} \circ \eta |^2_{\mathbb{H}}}{t-s}\right)| \varsigma |_{\mathbb{H}}^{\gamma} u^p(\varsigma , s) d\varsigma ds\\
 \noalign{\medskip}
	& \leq  C \sup_{t \in [0,T]}\| u(t) \|_{L^{\infty}(\mathbb{H}^N)}^p \int_{0}^{t} (t-s)^{-\frac{Q}{2}} \int_{\mathbb{H}^N}  \exp\left(  -c_\star \dfrac{| \varsigma^{-1} \circ \eta |^2_{\mathbb{H}}}{t-s}\right)| \varsigma |_{\mathbb{H}}^{\gamma} d\varsigma ds.
\end{array}
$$
Thus, from Lemma \ref{UpperSemigroup}, it follows that
\begin{align}\label{use-cal-frent2}
    \begin{aligned}
        [\mathcal{G}u(t)](\eta)  & \leq  C  K^p \int_{0}^{t} (t-s)^{\frac{\gamma}{2}}  ds \leq C_{5} K^p  T^{1 +\frac{\gamma}{2} },
    \end{aligned}
\end{align}
where $C_{5}>0$ is a constant. 
Therefore
\begin{align*}
	\begin{aligned}
		\left| 	[\mathcal{F}u(t)] (\eta)  \right| &\leq | [e^{ t\Delta_{\mathbb{H}} }u_0](\eta)| + | [\mathcal{G}u(t)](\eta) | \\
\noalign{\medskip}
   &\leq   \|u_0\|_{L^{\infty}(\mathbb{H}^N)}  \|h_t\|_{L^{1}(\mathbb{H}^N)} + C_5 K^p T^{1 + \frac{\gamma}{2} }\\
\noalign{\medskip}
   & \leq \|u_0\|_{L^{\infty}(\mathbb{H}^N)}  + C_5 K^p T^{1 + \frac{\gamma}{2} },
	\end{aligned}
\end{align*}
that is,
\begin{align}\label{Har1}
	\begin{aligned}
		\| 	\mathcal{F}u(t)  \|_ {L^{\infty}(\mathbb{H}^N)}\leq \|u_0\|_{L^{\infty}(\mathbb{H}^N)}  + C_5 K^p T^{1 + \frac{\gamma}{2} }.
	\end{aligned}
\end{align}

Similarly, for $u, v \in B_K$ , there is a constant $C_6>0$ such that
\begin{align} \label{Har2}
	\begin{aligned}
		\| 	\mathcal{F}u (t) - 	\mathcal{F}v (t) \|_ {L^{\infty}(\mathbb{H}^N)}\leq  C_6 K^{p-1} T^{1 +\frac{\gamma}{2} } 	\| 	u (t) - 	v (t) \|_ {L^{\infty}(\mathbb{H}^N)}.
	\end{aligned}
\end{align}

Now, we must choose $K > \|u_0\|_{L^{\infty}(\mathbb{H}^N)}  + 1$ and $T$ small enough so that $ C_6 K^{p-1} T^{1 +\frac{\gamma}{2} } <1$  and $C_5 K^p T^{1 + \frac{\gamma}{2} }< 1.$
Therefore, from \eqref{Har1} and \eqref{Har2}, it follows that $\mathcal{F}$ is a contraction in $B_K$ and possesses a unique fixed point $u$ in $B_K$.

The uniqueness follows from
$$\|u(t) - v(t)\|_ {L^{\infty}(\mathbb{H}^N)}\leq \sup_{s \in [0,T]} \left(  \|  u(s) \|_ {L^{\infty}(\mathbb{H}^N)}^{p-1} +  \|  v(s) \|_ {L^{\infty}(\mathbb{H}^N)}^{p-1} \right) \int_{0}^{t} (t-s)^{\gamma/2}\|u(s) - v(s)\|_ {L^{\infty}(\mathbb{H}^N)} ds $$
and the singular Gronwall inequality (see \cite[Proposition A.5.7]{breziscaz}) since ${\gamma}/{2}>-1.$ 

Since $T$ depends on $\| u_0\|_ {L^{\infty}(\mathbb{H}^N)}$, the proof of the blow-up alternative can be completed by standard arguments (see \cite[Remark 2.4.2 and Theorem 2.4.3]{breziscaz}). 
\cqd

\section{Global and blow-up results for the H\'{e}non equation}
\label{GloBl}
The following technical estimate will be used to prove the global existence of solutions. 
\begin{lemma} \label{Henon1} Let $\gamma \geq 0$, $u_0 \in L^\infty(\mathbb{H}^N) \cap L^1(\mathbb{H}^N)$, $u_0\geq 0$, $u_0 \neq 0$ such that $u_0(\eta) \leq (1+|\eta|_{\mathbb{H}})^{-Q}$. Then
$$\| \varphi \cdot e^{t \Delta_{\mathbb{H}}} u_0 \|_{L^\infty(\mathbb{H}^N)} \leq  C t^{ -\frac{Q}{2}+\frac{\gamma}{2(p-1)} },$$
for all $t>1$.
\end{lemma}
{\bf Proof.} Suppose first that $\gamma>0$. By estimate (\ref{estimativadokernel})
\begin{align}\label{eq-extra-expl}
\begin{aligned}
(1+|\eta|_{\mathbb{H}})^{\frac{\gamma}{p-1}} e^{t \Delta_{\mathbb{H}}} u_0(\eta) & \leq C_\star t^{-\frac{Q}{2}} 	(1+|\eta|_{\mathbb{H}})^{\frac{\gamma}{p-1}} \int_{\mathbb{H}^N} e^{-\frac{c_\star|\xi^{-1} \circ \eta|_{\mathbb{H}}^2 }{t}} u_0(\xi) d\xi \\
\noalign{\medskip}
& \leq C_\star t^{-\frac{Q}{2}} 	(1+|\eta|_{\mathbb{H}})^{\frac{\gamma}{p-1}} \left(  \int_{  |\xi|_{\mathbb{H}} \leq \frac{    | \eta |_{\mathbb{H}}    }{2}      } + \int_{  |\xi|_{\mathbb{H}} > \frac{    | \eta |_{\mathbb{H}}    }{2}      }  \right) e^{-\frac{c_\star|\xi^{-1} \circ \eta|_{\mathbb{H}}^2 }{t}} u_0(\xi) d\xi \\
\noalign{\medskip}
& = I+II.
\end{aligned}
\end{align}
From the triangle inequality, $|\xi^{-1}  \circ \eta |_{\mathbb{H}}  \geq  | \eta |_{\mathbb{H}} - |\xi |_{\mathbb{H}}$,  we have for $t>1$
\begin{align*}
I &= C_\star t^{-\frac{Q}{2}} 	(1+|\eta|_{\mathbb{H}})^{\frac{\gamma}{p-1}}   \int_{  |\xi|_{\mathbb{H}} \leq \frac{    | \eta |_{\mathbb{H}}    }{2}      }  e^{-\frac{c_\star|\xi^{-1} \circ \eta|_{\mathbb{H}}^2 }{t}} u_0(\xi) d\xi \\
\noalign{\medskip}
&\leq C_\star t^{-\frac{Q}{2} } 	(1+|\eta|_{\mathbb{H}})^{\frac{\gamma}{p-1}}  e^{-\frac{c_\star| \eta|_{\mathbb{H}}^2 }{4t}}   \int_{  |\xi|_{\mathbb{H}} \leq \frac{    | \eta |_{\mathbb{H}}    }{2}      }  u_0(\xi) d\xi \\
\noalign{\medskip}
&\leq C_\star t^{-\frac{Q}{2} } 	\left(\sqrt{t}+\sqrt{t} \, \frac{ |\eta|_{\mathbb{H}} }{\sqrt{t} } \right)^{\frac{\gamma}{p-1}}  e^{-\frac{c_\star| \eta|_{\mathbb{H}}^2 }{4t}}   \int_{  |\xi|_{\mathbb{H}} \leq \frac{    | \eta |_{\mathbb{H}}    }{2}      }  u_0(\xi) d\xi \\
\noalign{\medskip}
&\leq Ct^{-\frac{Q}{2}  + \frac{\gamma}{2(p-1)}} .\end{align*}	 
Now, using the assumption $ u_0(\xi) \leq (1+|\xi|_{\mathbb{H}})^{-Q} $ and arguing in a similar way to the derivation of the estimate \eqref{Upper1}, we obtain:
\begin{align*}
II &\leq C_\star t^{-\frac{Q}{2}} 	(1+|\eta|_{\mathbb{H}})^{\frac{\gamma}{p-1}}   \int_{  |\xi|_{\mathbb{H}} > \frac{    | \eta |_{\mathbb{H}}    }{2}      }  e^{-\frac{c_\star |\xi^{-1} \circ \eta|_{\mathbb{H}}^2 }{t}} u_0(\xi) d\xi \\
\noalign{\medskip}
&\leq 2^{\frac{\gamma}{p-1}} C_\star t^{-\frac{Q}{2}} 	 \int_{  |\xi|_{\mathbb{H}} > \frac{    | \eta |_{\mathbb{H}}    }{2}      }  e^{-\frac{c_\star|\xi^{-1} \circ \eta|_{\mathbb{H}}^2 }{t}} (1+|\xi|_{\mathbb{H}})^{-Q+\frac{\gamma}{p-1}}  d\xi \\
\noalign{\medskip}
&\leq 2^{\frac{\gamma}{p-1}} C_\star t^{-\frac{Q}{2}} 	 \int_{  |\xi|_{\mathbb{H}} > \frac{    | \eta |_{\mathbb{H}}    }{2}      }  e^{-\frac{c_\star |\xi^{-1} \circ \eta|_{\mathbb{H}}^2 }{t}} |\xi|_{\mathbb{H}}^{-Q+\frac{\gamma}{p-1}}  d\xi \\
\noalign{\medskip}
&\leq C t^{ \frac{-Q}{2}+\frac{\gamma}{2(p-1)} }.
\end{align*}	 
Therefore
$$\| \varphi \cdot e^{t \Delta_{\mathbb{H}}} u_0 \|_{L^\infty(\mathbb{H}^{N})} \leq  C t^{ -\frac{Q}{2}+\frac{\gamma}{2(p-1)} },$$
for $t>1$. 

The case $\gamma = 0$ follows immediately from the first line of the estimate \eqref{eq-extra-expl}, using the fact that $u_0 \in L^1(\mathbb{H}^N).$
\cqd   \\

The following lemma provides a sufficient condition for the global existence of solutions when $\gamma \geq 0$.
\begin{lemma}\label{GlobPos}
 Assume $p>1$ and $\gamma \geq 0$. If there is $ 0 \leq w_0 \in L^{\infty}_{\varphi}(\mathbb{H}^N),$ $w_0\neq 0$ such that 
 \begin{equation}\label{cond-sol-globais-heisen}
\int_{0}^{\infty}  \|  \varphi \cdot e^{s \Delta_{\mathbb{H}}} {w}_{0}  \|_{L^\infty(\mathbb{H}^{N})}^{p-1} ds < \infty,
\end{equation}
then there exists a global solution of \eqref{Eqgeral-1}.
\end{lemma}
{\bf Proof.} Let $u_0 = \lambda w_0$, where $0<\lambda \leq (1+\Lambda)^{-\frac{p}{p-1}}  $,  where  $\Lambda>0$ is given by
\begin{equation*}
\Lambda= \int_{0}^{\infty}  \|  \varphi \cdot e^{s \Delta_{\mathbb{H}}} {w}_{0}  \|_{L^\infty(\mathbb{H}^{N})}^{p-1} ds.
\end{equation*}
Consider the sequence $(u_n)_{n\geq 0}$ given by $u_{n=0}(t)=e^{t \Delta_{\mathbb{H}}} u_0$,
$$u_{n+1}(t)=e^{t \Delta_{\mathbb{H}}}u_0+\int_0^t e^{(t-s) \Delta_{\mathbb{H}}} \left[  |\cdot|_{\mathbb{H}}^{\gamma} u_n^p(s) \right] ds.$$
We claim that  
\begin{equation}\label{Cec.die-esph}
u_n(t)\leq (1+\Lambda)w(t),
\end{equation}
for $t>0$ and $n \geq 0$, where $w(t)= e^{t \Delta_{\mathbb{H}}}u_0 $. We argue by induction on $n$. Clearly, the estimate (\ref{Cec.die-esph}) holds  for $n=0$. Assuming that (\ref{Cec.die-esph}) is valid for $n$, we obtain
$$
\begin{array}{ll}
u_{n+1}(t)&= e^{t \Delta_{\mathbb{H}}} {u}_{0}+\int_{0}^{t} e^{(t-s) \Delta_{\mathbb{H}}} \left (|\cdot|_{\mathbb{H}}^{\gamma}u_n^{p}(s)\right )ds\\
\noalign{\medskip}
& \leq e^{t \Delta_{\mathbb{H}}} {u}_{0} + (1+\Lambda)^p\int_{0}^{t} e^{(t-s ) \Delta_{\mathbb{H}}} |\cdot|_{\mathbb{H}}^{\gamma}w^p(s) ds\\
\noalign{\medskip}
& \leq e^{t \Delta_{\mathbb{H}}} {u}_{0} + (1+\Lambda)^p\int_{0}^{t} e^{(t-s ) \Delta_{\mathbb{H}}} [|\cdot|_{\mathbb{H}}^{\frac{\gamma}{p-1}}w(s)]^{p-1} w(s) ds\\
\noalign{\medskip}
& \leq e^{t \Delta_{\mathbb{H}}} {u}_{0} + (1+\Lambda)^p\int_{0}^{t}  \| \varphi \cdot w(s)\|_{L^\infty(\mathbb{H}^{N})}^{p-1} 	e^{(t-s) \Delta_{\mathbb{H}}} w(s) ds\\
\noalign{\medskip}
& \leq e^{t \Delta_{\mathbb{H}}} {u}_{0} + (1+\Lambda)^p\int_{0}^{t}  \| \varphi \cdot w(s)\|_{L^\infty(\mathbb{H}^{N})}^{p-1} 	e^{(t-s) \Delta_{\mathbb{H}}} e^{s \Delta_{\mathbb{H}}} {u}_{0}  ds\\
\noalign{\medskip}
& \leq  e^{t \Delta_{\mathbb{H}}} {u}_{0}  + (1+\Lambda)^p  e^{ t \Delta_{\mathbb{H}}} {u}_{0} \int_0^t  \| \varphi \cdot w(s)\|_{L^\infty(\mathbb{H}^{N})}^{p-1} ds \\
\noalign{\medskip}
& \leq  w(t) + (1+\Lambda)^p  w(t) \lambda^{p-1} \int_{0}^{t}  \|  \varphi \cdot  e^{s \Delta_{\mathbb{H}}} {w}_{0}  \|_{L^\infty(\mathbb{H}^{N})}^{p-1} ds \\
\noalign{\medskip}
& \leq  \left[ 1+(1+\Lambda)^p \lambda^{p-1}  \Lambda \right] w(t)\\
\noalign{\medskip}
& = (1+\Lambda)w(t).
\end{array}
$$
Therefore, the claim (\ref{Cec.die-esph}) holds for $n+1$. 

Since $(u_n)_{n\geq 0}$ is clearly nondecreasing, we conclude, by the estimate (\ref{Cec.die-esph}), that there exists $u(t)=\lim_{n\to \infty} u_n(t)$. By the monotone convergence theorem, we conclude that $u$ is a global solution to problem \eqref{Eqgeral-1}.
\cqd
~~\\

For the proof of Theorems \ref{global-e-blowup-gamapositivo} (i) and \ref{globalsolgamaneg} (i), we will need the following lemma.

\begin{lemma}\label{lemma-tipo-fuj} Let $u$ be a solution of \eqref{Eqgeral-1} in $\mathbb{H}^N \times [0,T]$ with $0<T \leq \infty$.  

\begin{enumerate}[(i)]
    \item If $0\leq \gamma<Q(p-1)$, then there exists
a constant $M>0$ depending only on $p, \,Q$ and $ \gamma$ such that
\begin{equation}\label{eq-uso-caso-crit}
  s^{\frac{2+\gamma}{2}} \, \|e^{ s\Delta_{\mathbb{H}} }u_0\|_{L^\infty(\mathbb{H}^N)}^{p-1}  \leq M,
\end{equation}
for any $s \in [0,T)$. 

\item If $-2<\gamma<0$, then 
\begin{equation}\label{eq-22t11-b-neg}
\lim_{s \to \infty} \left\| \left[e^{ s\Delta_{\mathbb{H}} }u_0 \right]^{p-1}  \int_{0}^{s}\left [e^{(s-t)\Delta_{\mathbb{H}}}|\cdot|_{\mathbb{H}}^{-\frac{\gamma}{p-1}} \right]^{1-p}  dt \right\|_{L^\infty(\mathbb{H}^N)} < \infty.
\end{equation}
\end{enumerate}

\end{lemma}
{\bf Proof.}
Let $t\in (0,s)$, since 
$$ 
 \begin{array}{ll}
 	u(t) = e^{ t\Delta_{\mathbb{H}} }u_0 +   \int_{0}^{t} e^{(t-\sigma)\Delta_{\mathbb{H}}}| \cdot |_{\mathbb{H}}^{\gamma} u^p(\sigma) d\sigma ,
 \end{array}
 $$  
using Fubini's Theorem, we can apply $e^{(s-t)\Delta_{\mathbb{H}}}$ to both sides and we get
$$
 \begin{array}{ll}
 	e^{(s-t)\Delta_{\mathbb{H}}}u(t)= e^{ s\Delta_{\mathbb{H}} }u_0 + \int_{0}^{t} e^{(s-\sigma)\Delta_{\mathbb{H}}}| \cdot |_{\mathbb{H}}^{\gamma} u^p(\sigma) d\sigma. \\
\end{array}
$$


\noindent{\bf (i)} $0\leq \gamma<Q(p-1)$. Similarly to the derivation of the estimate \eqref{Upper1}, we can infer that
$$ \int_{\mathbb{H}^N} h_{s-\sigma}(\varsigma^{-1} \circ \eta) |\varsigma|_{\mathbb{H}}^{- \frac{\gamma}{p-1}} d \varsigma \in (0, \infty).$$
 According to the reverse Hölder inequality with $f(\varsigma) = \left[ h_{s-\sigma}(\varsigma^{-1} \circ \eta) u(\sigma, \varsigma) \right]^p  $ and $g(\varsigma) = \left[ h_{s-\sigma}(\varsigma^{-1} \circ \eta) |\varsigma|_{\mathbb{H}}^{-\frac{\gamma}{p-1}}  \right]^{1-p}$
(Lemma \ref{invholder}), we can obtain
$$e^{(s-t)\Delta_{\mathbb{H}}} u(t) \geq  e^{ s\Delta_{\mathbb{H}} }u_0 + \int_{0}^{t} \left [e^{(s-\sigma)\Delta_{\mathbb{H}}}| \cdot |_{\mathbb{H}}^{-\frac{\gamma}{p-1}} \right]^{1-p} \left [e^{(s-\sigma)\Delta_{\mathbb{H}}}u(\sigma) \right ]^{p} d\sigma.$$
Set  
 $$ \Phi(t) := \Phi(t, \cdot) =  e^{ s\Delta_{\mathbb{H}} }u_0 + \int_{0}^{t} \left [e^{(s-\sigma)\Delta_{\mathbb{H}}}| \cdot |_{\mathbb{H}}^{-\frac{\gamma}{p-1}} \right]^{1-p} \left [e^{(s-\sigma)\Delta_{\mathbb{H}}}u(\sigma) \right ]^{p} d\sigma.$$
We have that $\Phi(t, \cdot) \leq e^{(s-t)\Delta_{\mathbb{H}}} u(t) \leq \|h_{s-t}\|_{L^1(\mathbb{H}^N)} \| u(t)\|_{L^\infty(\mathbb{H}^N)} \leq \|\varphi \cdot  u(t)\|_{L^\infty(\mathbb{H}^N)}<\infty$ and it is absolutely continuous on $[0,s].$  Hence
 \begin{align*}
 	\Phi'(t) &=  \left [e^{(s-t)\Delta_{\mathbb{H}}}|\cdot|_{\mathbb{H}}^{-\frac{\gamma}{p-1}} \right]^{1-p} \left [e^{(s-t)\Delta_{\mathbb{H}}}u(t) \right ]^{p}\\ 
  \noalign{\medskip}
 	&\geq  \left [e^{(s-t)\Delta_{\mathbb{H}}}|\cdot|_{\mathbb{H}}^{-\frac{\gamma}{p-1}} \right]^{1-p}  [\Phi(t)]^{p},
 \end{align*}
 and so
 \begin{align*}
 	(\Phi^{1-p})'(t) &\leq -(p-1)   \left [e^{(s-t)\Delta_{\mathbb{H}}}|\cdot|_{\mathbb{H}}^{-\frac{\gamma}{p-1}} \right]^{1-p} .
 \end{align*}
 Thus, integrating on $[0,s]$ and using the Fundamental Theorem of Calculus, we obtain
$$
 \begin{array}{ll}
 	(\Phi^{1-p})(0) &\geq 	(\Phi^{1-p})(s)  +  (p-1) \int_{0}^{s}\left [e^{(s-t)\Delta_{\mathbb{H}}}|\cdot|_{\mathbb{H}}^{-\frac{\gamma}{p-1}} \right]^{1-p}  dt\\
  \noalign{\medskip}
 	&\geq   (p-1) \int_{0}^{s}\left [e^{(s-t)\Delta_{\mathbb{H}}}|\cdot|_{\mathbb{H}}^{-\frac{\gamma}{p-1}} \right]^{1-p}  dt.
 \end{array}
$$ 
 Therefore
\begin{equation}\label{cor.onc} 
 [e^{ s\Delta_{\mathbb{H}} }u_0]^{p-1}  \int_{0}^{s}\left [e^{(s-t)\Delta_{\mathbb{H}}}|\cdot|_{\mathbb{H}}^{-\frac{\gamma}{p-1}} \right]^{1-p}  dt \leq \dfrac{1}{p-1}
 \end{equation}
 for any $0 < t < s$. From \eqref{Upper1}, we obtain 
 $$
\begin{array}{ll}
e^{(s-t)\Delta_{\mathbb{H}}}|\cdot|_{\mathbb{H}}^{-\frac{\gamma}{p-1}}&\leq C (s-t)^{-\frac{\gamma}{2(p-1)}},
\end{array} 
$$
then replacing this estimate into (\ref{cor.onc}) we obtain
\begin{equation*}
    \|e^{ s\Delta_{\mathbb{H}} }u_0\|_{L^\infty(\mathbb{H}^N)}^{p-1} \int_0^s (s-t)^{\frac{\gamma}{2}} dt \leq C<\infty.
\end{equation*}

\noindent {\bf (ii)} $-2<\gamma<0$. As in the previous case, we show first that $$ \int_{\mathbb{H}^N} h_{s-\sigma}(\varsigma^{-1} \circ \eta) |\varsigma|_{\mathbb{H}}^{- \frac{\gamma}{p-1}} d \varsigma \in (0, \infty).$$
By the symmetry and triangle inequality of $|\cdot|_{\mathbb{H}}$, we have $|\varsigma|_{\mathbb{H}}^{- \frac{\gamma}{p-1}} \leq C|\varsigma^{-1} \circ \eta|_{\mathbb{H}}^{- \frac{\gamma}{p-1}} + C|\eta|_{\mathbb{H}}^{- \frac{\gamma}{p-1}}$, where $C>0$ depends only on $\gamma, \, p$. Moreover, by \cite[Proposition 5.4.4.]{BONFIGLIOLI} concerning the integration of $|\cdot|_{\mathbb{H}}$-radially-symmetric functions, we have that
$$ \int_{B_{\mathbb{H}}(0,R)} \exp \left( -c_\star \frac{|\varsigma|^2_{\mathbb{H}}}{s-\sigma} \right)|\varsigma|_{\mathbb{H}}^{- \frac{\gamma}{p-1}} d \varsigma  = Q |B_{\mathbb{H}}(0,1)|  \int_{0}^R x^{Q- \frac{\gamma}{p-1}-1 } \exp \left( -c_\star \frac{x^2}{s-\sigma} \right) dx. $$
Thus, we can infer that
\begin{align*}
    \int_{\mathbb{H}^N} \exp \left( -c_\star \frac{|\varsigma^{-1} \circ \eta|^2_{\mathbb{H}}}{s-\sigma} \right)|\varsigma^{-1} \circ \eta|_{\mathbb{H}}^{- \frac{\gamma}{p-1}} d \varsigma  &= Q |B_{\mathbb{H}}(0,1)|  \int_{0}^\infty x^{Q- \frac{\gamma}{p-1}-1 } \exp \left( -c_\star \frac{x^2}{s-\sigma} \right) dx\\
    &\leq C (s - \sigma)^{\frac{Q}{2}- \frac{\gamma}{2(p-1)} }  .
\end{align*}
Therefore, it follows from this that
 \begin{align*}
     \int_{\mathbb{H}^N} h_{s-\sigma}(\varsigma^{-1} \circ \eta) |\varsigma|_{\mathbb{H}}^{- \frac{\gamma}{p-1}} d \varsigma &\leq C(s-\sigma)^{ - \frac{Q}{2}}   \int_{\mathbb{H}^N} \exp \left( -c_\star \frac{|\varsigma^{-1} \circ \eta|^2_{\mathbb{H}}}{s-\sigma} \right)|\varsigma|_{\mathbb{H}}^{- \frac{\gamma}{p-1}} d \varsigma\\
     &\leq C(s-\sigma)^{ - \frac{Q}{2}}   \int_{\mathbb{H}^N} \exp \left( -c_\star \frac{|\varsigma^{-1} \circ \eta|^2_{\mathbb{H}}}{s-\sigma} \right)|\varsigma^{-1} \circ \eta|_{\mathbb{H}}^{- \frac{\gamma}{p-1}} d \varsigma \\
     & \qquad + C|\eta|_{\mathbb{H}}^{- \frac{\gamma}{p-1}}(s-\sigma)^{ - \frac{Q}{2}}   \int_{\mathbb{H}^N} \exp \left( -c_{\star} \frac{|\varsigma^{-1} \circ \eta|^2_{\mathbb{H}}}{s-\sigma} \right) d \varsigma\\
     & \leq C \left [(s - \sigma)^{- \frac{\gamma}{2(p-1)} } + |\eta|_{\mathbb{H}}^{- \frac{\gamma}{p-1}} \right ]< \infty.
 \end{align*}

Proceeding as in the previous case we can obtain (\ref{cor.onc}) and the result follows. 

\cqd
~~\\

{\bf Proof of Theorem \ref{global-e-blowup-gamapositivo} } \textbf{(i) Case $p<p_c$. } 
We argue by contradiction. 
Consider any $u_0\geq 0$ such that  $u_0\geq \varepsilon>0$ in some ball $B_{\mathbb{H}}(\eta_0,r)$, $r>0$. Replace $u_0$ with a function $\tilde{u}_0 \in C_c(\mathbb{H}^N)$, where $C_c(\mathbb{H}^N)$ denotes the set of continuous functions with compact support in $\mathbb{H}^N$, satisfying $u_0\geq \tilde{u}_0 \geq \varepsilon>0$ in $B_{\mathbb{H}}(\eta_0,r)$. Then, by estimate \eqref{estimativadokernel} we have  
$$
\begin{array}{l}
s^{\frac{2+\gamma}{2}}\, \|e^{ s\Delta_{\mathbb{H}} }u_0\|_{L^\infty(\mathbb{H}^N)}^{p-1}  \\
\noalign{\medskip}
\geq s^{\frac{2+\gamma}{2}}\,\|e^{ s\Delta_{\mathbb{H}} }\tilde{u}_0\|_{L^\infty(\mathbb{H}^N)}^{p-1} \\
\noalign{\medskip}
\geq s^{\frac{2+\gamma}{2}}\, \left[ e^{ s\Delta_{\mathbb{H}} }\tilde{u}_0 (0) \right]^{p-1}  \\ 
\noalign{\medskip}
\geq  c_\star^{p-1} s^{-\frac{Q}{2}(p-1) + \frac{2+\gamma}{2}} \left( \int_{\mathbb{H}^N} e^{-\frac{C_{\star}|\xi|_{\mathbb{H}}^2 }{s}} \tilde{u}_0(\xi) d\xi  \right)^{p-1}\\
\noalign{\medskip}
\geq  c_\star^{p-1} s^{-\frac{Q}{2}(p-1) + \frac{2+\gamma}{2}}  \left(\int_{B_{\mathbb{H}}(\eta_0,r)} e^{-{C_\star|\xi|_{\mathbb{H}}^2 }} \tilde{u}_0(\xi) d\xi \right)^{p-1} \\
\noalign{\medskip}
\geq C s^{-\frac{Q}{2}(p-1) + \frac{2+\gamma}{2}} \to +\infty
\end{array}
$$	
when $s\to +\infty$ since $1<p<1+(2+\gamma)/Q.$ This contradicts (\ref{eq-uso-caso-crit}) in Lemma \ref{lemma-tipo-fuj} $(i)$.  \\

\noindent \textbf{ Critical case $p=p_c$. } We will argue by contradiction and suppose that there exists $u\in L^{\infty}([0,\infty), L^{\infty}_{\varphi}(\mathbb{H}^N))$ a global-in-time solution (i.e., $T_{\max} = \infty$) of problem \eqref{Eqgeral-1} with the initial condition $u_0 \in L^{\infty}_{\varphi}(\mathbb{H}^N)$. We can also assume that $u_0>0$. Due to the Gaussian estimates \eqref{estimativadokernel} and the triangle inequality for $|\cdot|_{\mathbb{H}}$, we have  
\begin{align}\label{eqfc3-A1}
\begin{aligned}
        e^{t \Delta_{\mathbb{H}}}u_{0}(\eta)
        &\geq 
    \int_{\mathbb{H}^N} h_t(\varsigma^{-1} \circ \eta) u_0(\varsigma) d \varsigma\\
    \noalign{\medskip}
   & \geq c_\star t^{-\frac{Q}{2}} \exp \left( -\frac{2C_{\star} |\eta|_{\mathbb{H}}^2}{t}\right)\int_{| \varsigma|_{\mathbb{H}} \leq \sqrt{t}}  \exp \left( -\frac{2C_{\star} |\varsigma|_{\mathbb{H}}^2}{t}\right) u_{0}(\varsigma)  d\varsigma\\
   \noalign{\medskip}
   & \geq c _\star \exp(-2C_{\star}) \exp \left( -\frac{2C_{\star} |\eta|_{\mathbb{H}}^2}{t} \right) t^{-\frac{Q}{2}} \int_{| \varsigma|_{\mathbb{H}} \leq 1} u_{0}(\varsigma)  d\varsigma\\
   &=Ct^{-\frac{Q}{2}} \exp \left( -\frac{2C_{\star} |\eta|_{\mathbb{H}}^2}{t} \right) ,
\end{aligned}
\end{align}
for $t \geq  1$, where $C = c_\star \exp(-2C_{\star}) \int_{| \varsigma|_{\mathbb{H}} \leq 1} u_{0}(\varsigma)  d\varsigma>0 $. Since $ t + 1 - s \leq t $ and $ s \leq t + 1 - s $ for $ 1 \leq s \leq t/2$, by  the equality (\ref{defmildsol}),  estimate (\ref{eqfc3-A1}) and \cite[Proposition 5.4.4.]{BONFIGLIOLI} we have
\begin{align}\label{princ-arFC}
    \begin{aligned}
&\int_{|\eta|_{\mathbb{H}}\leq \sqrt{t+1}}u(\eta ,t+1)d\eta \\
\noalign{\medskip}
& \geq \int_{|\eta|_{\mathbb{H}}\leq \sqrt{t}} \int_{0}^{t+1} e^{(t+1 - s)\Delta_{\mathbb{H}}}  \left( |\cdot|_{ \mathbb{H} }^{\gamma}u^{p}(s) \right) (\eta) ds \,d\eta \\
\noalign{\medskip}
&\geq  c_*\int_{|\eta|_{\mathbb{H}}\leq \sqrt{t}} \int_{1}^{t/2}  (t+1 - s)^{- \frac{Q}{2}} \int_{\mathbb{H}^N} e^{-\frac{C_* |\varsigma^{-1}\circ \eta|^2_{\mathbb{H}}}{t+1 - s}} |\varsigma|_{ \mathbb{H} }^{\gamma} u^p(\varsigma,s)  d \varsigma \, d s \, d\eta \\
&\geq  c_* C^p\int_{|\eta|_{\mathbb{H}}\leq \sqrt{t}} \int_{1}^{t/2}  (t+1 - s)^{- \frac{Q}{2}} s^{-\frac{Q}{2}p}\int_{\mathbb{H}^N} e^{-\frac{C_* |\varsigma^{-1}\circ \eta|^2_{\mathbb{H}}}{t+1 - s}} |\varsigma|_{ \mathbb{H} }^{\gamma} e^{-\frac{2pC_*|\varsigma|^2_{\mathbb{H}}}{s}}   d \varsigma \, d s \, d\eta \\
 &\geq  c_* C^p \int_{1}^{t/2}  (t+1 - s)^{- \frac{Q}{2}} s^{-\frac{Q}{2}p} \int_{|\eta|_{\mathbb{H}}\leq \sqrt{t}} e^{-\frac{2C_* |\eta|^2_{\mathbb{H}}}{t+1 - s}} \int_{\mathbb{H}^N} e^{-\frac{2C_* |\varsigma|^2_{\mathbb{H}}}{t+1 - s}} |\varsigma|_{ \mathbb{H} }^{\gamma} e^{-\frac{2pC_*|\varsigma|^2_{\mathbb{H}}}{s}}   d \varsigma \,  d\eta \, d s  \\
  &\geq  c_* C^p \int_{1}^{t/2}  (t+1 - s)^{- \frac{Q}{2}} s^{-\frac{Q}{2}p} \left( \int_{|\eta|_{\mathbb{H}}\leq \sqrt{t+1-s}} e^{-\frac{2C_* |\eta|^2_{\mathbb{H}}}{t+1 - s}} \,d\eta \right) \int_{|\varsigma|_{\mathbb{H}}\leq \sqrt{t+1-s}} e^{-\frac{2C_* |\varsigma|^2_{\mathbb{H}}}{t+1 - s}} |\varsigma|_{ \mathbb{H} }^{\gamma} e^{-\frac{2pC_*|\varsigma|^2_{\mathbb{H}}}{s}}   d \varsigma \,   d s  \\
       &\geq  c_* C^p C_1 \int_{1}^{t/2}  (t+1 - s)^{- \frac{Q}{2}} \, (t+1 - s)^{ \frac{Q}{2}}\, s^{-\frac{Q}{2}p} \int_{|\varsigma|_{\mathbb{H}}\leq \sqrt{t+1-s}} |\varsigma|_{ \mathbb{H} }^{\gamma} e^{-\frac{2pC_*|\varsigma|^2_{\mathbb{H}}}{s}}   d \varsigma \,   d s  \\   
         &\geq  c_* C^p C_1 \int_{1}^{t/2}  s^{-\frac{Q}{2}p} \int_{|\varsigma|_{\mathbb{H}}\leq \sqrt{s}} |\varsigma|_{ \mathbb{H} }^{\gamma} e^{-\frac{2pC_*|\varsigma|^2_{\mathbb{H}}}{s}}   d \varsigma \,   d s  \\ 
           &\geq  c_* C^p C_1 Q |B_{\mathbb{H}}(0,1)| \int_{1}^{t/2}  s^{-\frac{Q}{2}p} s^{\frac{Q}{2} +\frac{\gamma}{2} } \left( \int_0^1 x^{Q+\gamma-1} e^{-2pC_* \,x^2}   d x \right) \,   d s  \\ 
           & =   c_* C^p C_1C_2 Q |B_{\mathbb{H}}(0,1)| \int_{1}^{t/2}  s^{-\frac{Q}{2}(p-1) +\frac{\gamma}{2}  } \,   d s  \\
             & =   c_* C^p C_1C_2 Q |B_{\mathbb{H}}(0,1)| \int_{1}^{t/2}  s^{-1 }   \,   d s  \\
           &\geq C_3 \ln \left( \frac{t}{2} \right) \to \infty,
\end{aligned}
	\end{align}
as $t \to \infty$, where $C_1 = e^{-4C_*}$, $C_2 = \textstyle \int_0^1 x^{Q+\gamma-1} e^{-2pC_* \,x^2}   d x $ and $C_3 = c_* C^p C_1C_2 Q |B_{\mathbb{H}}(0,1)|$.

 Let $L$ be a positive constant that will be chosen later to be sufficiently large. By \eqref{princ-arFC} we can find $T_0 > 2$ sufficiently large such that the function $U_{0}$ defined by $U_{0} (\eta):= u(\eta, T_0 ) \ $ satisfies
$$\int_{|\eta|_{\mathbb{H}}\leq \sqrt{T_0}} U_{0}(\eta)d\eta \geq L.$$ 
It follows from this and from the same argument used in \eqref{eqfc3-A1} that
$$[e^{t\Delta_{\mathbb{H}}} U_{0}](\xi) \geq C_1 t^{-\frac{Q}{2}}\int_{|\eta|_{\mathbb{H}}\leq \sqrt{T}} U_{0}(\eta)d\eta \geq C_1L\,t^{-\frac{Q}{2}}$$
for all $t \geq T_0$ and $|\xi|_{\mathbb{H}} \leq \sqrt{t}$. This implies that 
\begin{equation}\label{pcontr1}
t^{\frac{Q}{2}}\| e^{t\Delta_{\mathbb{H}}} U_{0} \|_{  L^{\infty}(\mathbb{H}^N) }  \geq   C_1L.
\end{equation}
On the other hand, let $v$ be a solution of \eqref{Eqgeral-1} with initial data $U_0(\cdot) = u(\cdot, T_0)$ (defined above). Then, by uniqueness, since $u$ is a global solution of \eqref{Eqgeral-1}, $v$  is also a global solution of \eqref{Eqgeral-1}. Therefore, we can apply Lemma \ref{lemma-tipo-fuj} $(i)$ to the solution $v$, and since $Q/2 = (2+\gamma)/(2(p-1))$, we have  
\begin{equation}\label{pcontr2}
    M \geq t^{ \frac{2+\gamma}{2(p-1)} } \|e^{ t\Delta_{\mathbb{H}} }U_0\|_{L^\infty(\mathbb{H}^N)} = t^{ \frac{Q}{2} } \|e^{ t\Delta_{\mathbb{H}} }U_0\|_{L^\infty(\mathbb{H}^N)} 
\end{equation}
for any  $t > 0$. Recalling that $M$ depends on $p,Q$ and $\gamma$. Therefore, it follows from \eqref{pcontr1} and \eqref{pcontr2} that $M \geq C_1 L.$ This contradicts the arbitrariness of $L$. Thus, the problem \eqref{Eqgeral-1} does not have global nonnegative nontrivial solutions. \\
~~\\
\medskip
\textbf{(ii)}  Observing the condition \eqref{cond-sol-globais-heisen} of Lemma \ref{GlobPos}, we can see that if $\int_{0}^{\infty} \| \varphi \cdot e^{s \Delta_{\mathbb{H}}} {w}_{0} \|_{L^\infty(\mathbb{H}^N)}^{p-1} ds< \infty$, then it is possible to construct a global solution of \eqref{Eqgeral-1}. For $0<t<1$, the estimate \eqref{use-more-n=one} shows that $\int_{0}^{1} \| \varphi \cdot e^{s \Delta_{\mathbb{H}}} {w}_{0} \|_{L^\infty(\mathbb{H}^N)}^{p-1} ds< \infty$. Also, from Lemma \ref{Henon1}, for $t>1$ we have the following estimate
$$\| \varphi \cdot  e^{s \Delta_{\mathbb{H}}} w_0 \|_{L^\infty(\mathbb{H}^N)} \leq C s^{ -\frac{Q}{2}+\frac{\gamma}{2(p -1)} }$$
for $w_0 \in C_c(\mathbb{H}^N)$, $w_0(\eta)\geq 0$, $w_0 \neq 0$ such that $w_0(\eta) \leq [1+|\eta |_{\mathbb{H}}]^{-Q}$ (note that $ w_0 $ satisfying these conditions, belongs to $ L^{\infty}_{\varphi}(\mathbb{H}^N)$). Therefore we conclude that the integral in \eqref{cond-sol-globais-heisen} is finite, if
$$ \int_{1}^{\infty} s^{ -\frac{Q(p-1)}{2}+\frac{\gamma}{2} } ds< \infty.$$
But this occurs if $-{Q(p-1)}/{2}+{\gamma}/{2}<-1$, that is, if $ p> 1+ {(2+\gamma) }/{Q} .$
\cqd
~~\\

\section{Global and blow-up results for the Hardy equation}
\label{GloB2}

In order to prove the global existence of solutions, we will require the following lemma.
\begin{lemma} \label{GlobNeg}
Assume $p>1$, $\gamma<0$, and $q>1$ such that $1/q<1+\gamma /Q$. If there exists $w_0 \in L^{\infty}(\mathbb{H}^N)$, $w_0 \geq 0$, $w_0\neq 0$ such that
\begin{equation}\label{usaas}
     \Lambda= \sup_{t>0}\int_{0}^{t}  \|  e^{\sigma \Delta_{\mathbb{H}} }w_{0}  \|_{L^\infty(\mathbb{H}^N)}^{\frac{p-1}{q}} (t-\sigma)^{\frac{\gamma}{2}} d\sigma<\infty,
\end{equation}
then there exists a global solution of \eqref{Eqgeral-1}.
\end{lemma}
{\bf Proof.} Let $q'>1$ such that $ {1}/{q} + {1}/{q'}=1$, and
\begin{equation}\label{usaas}
     \Lambda= \sup_{t>0}\int_{0}^{t}  \|  e^{\sigma \Delta_{\mathbb{H}} }w_{0}  \|_{L^\infty(\mathbb{H}^N)}^{\frac{p-1}{q}} (t-\sigma)^{\frac{\gamma}{2}} d\sigma<\infty.
\end{equation}

By Lemma \ref{UpperSemigroup}, there is a positive constant $C_0>0$ so that
\begin{equation} \label{Gamma1}
e^{(t-\sigma) \Delta_{\mathbb{H}}} |\cdot|^{\gamma q'}_{\mathbb{H}} \leq C_0 (t-\sigma)^{ \frac{\gamma q'}{2}},
\end{equation}
for $0< \sigma < t.$

Consider $u_0 = \lambda w_0$, where $0<\lambda < \min \left\{ C^{-\frac{q}{q'(p-1)}}_{0}(1+\Lambda)^{-\frac{pq}{p-1}} , \|w_0\|_{L^\infty(\mathbb{H}^N)}^{-1} \right\}$. We have
$$ e^{t \Delta_{\mathbb{H}} }u_0\leq \lambda   \|h_t\|_{L^{1}(\mathbb{H}^N)} \|w_0\|_{L^\infty(\mathbb{H}^N)} < 1,$$
and therefore, $e^{t \Delta_{\mathbb{H}} }u_0 \leq (e^{t \Delta_{\mathbb{H}} }u_0)^{\frac{1}{q}}$. 

Now, we define the sequence $(u_n)_{n\geq 0}$ given by $u_{n=0}(t)=e^{t \Delta_{\mathbb{H}} }u_0$,
\begin{equation}\label{Sh.uno}
u_{n+1}(t)=e^{t \Delta_{\mathbb{H}} }u_0 +\int_0^te^{(t - \sigma) \Delta_{\mathbb{H}} } \left[  |\cdot|^{\gamma}_{\mathbb{H}} u_n^p(\sigma) \right] d\sigma,
\end{equation}
and $w(t)= (e^{t \Delta_{\mathbb{H}} }u_0)^{\frac{1}{q}}$. We claim that  
\begin{equation}\label{Cec.die-11}
u_n(t)\leq (1+\Lambda)w(t).
\end{equation}
We argue by induction on $n$. Clearly, the estimate (\ref{Cec.die-11}) holds  for $n=0$. Assume that it holds for $n$. By H\"older's inequality and \eqref{Gamma1}, we have
$$
\begin{array}{ll}
&u_{n+1}(t)\\
\noalign{\medskip}
&= e^{t \Delta_{\mathbb{H}} }u_0 +\int_0^t e^{(t - \sigma) \Delta_{\mathbb{H}} } \left[  |\cdot|^{\gamma}_{\mathbb{H}} u_n^p(\sigma) \right] d\sigma\\
\noalign{\medskip}
& \leq e^{t \Delta_{\mathbb{H}} }u_0 + (1+\Lambda)^p\int_{0}^{t} e^{(t - \sigma) \Delta_{\mathbb{H}} } (|\cdot|^{\gamma}_{\mathbb{H}}w^p(\sigma)) d\sigma\\
\noalign{\medskip}
& \leq e^{t \Delta_{\mathbb{H}} }u_0 + (1+\Lambda)^p\int_{0}^{t} \|w(\sigma)\|_{L^\infty(\mathbb{H}^N)}^{p-1} 	e^{(t - \sigma) \Delta_{\mathbb{H}} } |\cdot |^{\gamma}_{\mathbb{H}} (e^{\sigma \Delta_{\mathbb{H}} }u_0)^{\frac{1}{q}} d\sigma \\
\noalign{\medskip}
&\leq e^{t \Delta_{\mathbb{H}} }u_0 + (1+\Lambda)^p \int_0^t  \|w(\sigma)\|^{p-1}_{L^\infty(\mathbb{H}^N)} \left[ e^{(t - \sigma) \Delta_{\mathbb{H}} } |\cdot|^{\gamma q'}_{\mathbb{H}}\right]^{\frac{1}{q'}} \left[e^{(t - \sigma) \Delta_{\mathbb{H}} } (e^{\sigma \Delta_{\mathbb{H}} }u_0) \right]^{\frac{1}{q}}d\sigma
\\
\noalign{\medskip}
& \leq  e^{t \Delta_{\mathbb{H}} }u_0  + (1+\Lambda)^p \, [e^{t \Delta_{\mathbb{H}} }u_0 ]^{\frac{1}{q}} \int_0^t  \|w(\sigma)\|_{L^\infty(\mathbb{H}^N)}^{p-1} \left[e^{(t - \sigma) \Delta_{\mathbb{H}} } |\cdot|_{\mathbb{H}}^{\gamma q'} \right]^{\frac{1}{q'}} d\sigma \\
\noalign{\medskip}
& \leq  e^{t \Delta_{\mathbb{H}} }u_0 +  C^{\frac{1}{q'}}_{0} \lambda^{\frac{p-1}{q}} \, (1+\Lambda)^p \, [e^{t \Delta_{\mathbb{H}} }u_0 ]^{\frac{1}{q}} \,  \int_{0}^{t}  \| e^{\sigma \Delta_{\mathbb{H}} }w_{0}  \|_{L^\infty(\mathbb{H}^N)}^{\frac{p-1}{q}} (t-\sigma)^{\frac{\gamma}{2}}d\sigma \\
\noalign{\medskip}
& \leq  e^{t \Delta_{\mathbb{H}} }u_0 +  C^{\frac{1}{q'}}_{0} \lambda^{\frac{p-1}{q}} \, (1+\Lambda)^p \, \Lambda \, [e^{t \Delta_{\mathbb{H}} }u_0 ]^{\frac{1}{q}}  \\
\noalign{\medskip}
& = (1+\Lambda)w(t).
\end{array}
$$
Therefore, the claim (\ref{Cec.die-11}) holds for $n+1$. 
Since $(u_n)_{n\geq 0}$ is clearly nondecreasing we conclude, by the estimate (\ref{Cec.die-11}), that there exists $u(t)=\lim_{n\to \infty} u_n(t)$ and $u$ is a global solution to the problem (\ref{Eqgeral-1}).  \cqd \\ 


\medskip \noindent
{\bf Proof of Theorem \ref{globalsolgamaneg} }\textbf{ (i)}.  We argue by contradiction and suppose that there is a global solution $u$ of \eqref{Eqgeral-1} with the initial datum $u_0 \neq 0$ and $u_0\geq 0$. 

Let $\varepsilon>0$ be such that $u_0\geq \varepsilon>0$ in some ball $B_{\mathbb{H}}(\eta_0,r) \, (r>0)$ and $| \eta |_{\mathbb{H}} = ( |x|^2 + |y|^2 + |\tau|  )^{1/2} \, \, (\mbox{where } \eta =(x,y,\tau))$. Then, based on the estimate \eqref{estimativadokernel} we have
\begin{equation} \label{N1}
\begin{array}{l}
\left\| [e^{ s\Delta_{\mathbb{H}} }u_0]^{p-1}  \int_{0}^{s}\left [e^{(s-t)\Delta_{\mathbb{H}}}|\cdot|_{\mathbb{H}}^{-\frac{\gamma}{p-1}} \right]^{1-p}  dt \right\|_{L^\infty(\mathbb{H}^N)}  \\
\noalign{\medskip}
\geq \left ( [e^{ s\Delta_{\mathbb{H}} }u_0]^{p-1}  \int_{0}^{s}\left [e^{(s-t)\Delta_{\mathbb{H}}}|\cdot|_{\mathbb{H}}^{-\frac{\gamma}{p-1}} \right]^{1-p}  dt \right )(0)  \\
\noalign{\medskip}
\geq [e^{ s\Delta_{\mathbb{H}} }u_0 (0)]^{p-1}  \int_{0}^{s}\left [e^{(s-t)\Delta_{\mathbb{H}}}|\cdot|_{\mathbb{H}}^{-\frac{\gamma}{p-1}} \right]^{1-p} (0) dt\\
\noalign{\medskip}
\geq  c_{\star}^{p-1}s^{-\frac{Q}{2}(p-1)} 
\left(\int_{B_{\mathbb{H}}(\eta_0,r)} e^{-{C_\star \frac{ |\xi|_{\mathbb{H}}^2  }{s} }} u_0(\xi) d\xi \right)^{p-1} \cdot  \\
 \hskip20pt \int_{0}^{s} (s-t)^{-\frac{Q}{2}(1-p)}  \left[   \underbrace{ \int_{\mathbb{H}^N} e^{-{c_{\star} \frac{ |\eta |_{\mathbb{H}}^2  }{s-t} } } |\eta|^{ - \frac{\gamma}{p-1}}_{\mathbb{H}}  d\eta }_{J}  \right]^{1-p}  dt  \\
\noalign{\medskip}
\end{array}
\end{equation}
for $s>1$. Now, we will estimate $J$ considering the changes of variables $\tau' = (s- t) \tau
$, $y' = (s-t)^{\frac{1}{2}}y$, and $x' = (s-t)^{\frac{1}{2}}x$. Thus,
\begin{equation*}
\begin{array}{lll}
J & = &  \int_{\mathbb{R}^N}   \exp\left(  -c_\star \dfrac{  |x' |^2 }{s-t}\right)  \\
& & \times  \left [ \int_{\mathbb{R}^N}  \exp\left(  -c_\star \dfrac{ |y'|^2 }{s-t}\right)  \left [
\int_{\mathbb{R}}   \exp \left(  \dfrac{-c_\star | \tau' | }{s-t}  \right) ( |x|^2 + |y|^2 + |\tau|  )^{-\frac{\gamma}{2(p-1)}}     d\tau' \right  ] dy' \right ] dx' \\
\noalign{\medskip}
& =& (s-t)^{\frac{Q}{2}}  (s-t)^{\frac{-\gamma}{2(p-1)}}  \int_{\mathbb{R}^N} \exp\left(  -c_\star   |x |^2 \right)  \\
& & \times  \left [ \int_{\mathbb{R}^N} \exp\left(  -c_\star  |y|^2 \right)  \left [
\int_{\mathbb{R}} \exp \left( -c_\star | \tau | \right) ( [ |x|^2 + |y|^2 + |\tau| ] )^{-\frac{\gamma}{2(p-1)}}     d\tau \right  ] dy \right ] dx \\
\noalign{\medskip}
&=& (s-t)^{\frac{Q}{2}-\frac{\gamma}{2(p-1)}}   \int_{\mathbb{H}^N} e^{-c_\star  |\eta |_{\mathbb{H}}^2   }  |\eta |^{ - \frac{\gamma}{p-1}}_{\mathbb{H}}  d\eta.
\end{array}
\end{equation*}
Substituting this into estimate \eqref{N1}, we get
\begin{equation*} 
\begin{array}{l}
\left\| [e^{ s\Delta_{\mathbb{H}} }u_0]^{p-1}  \int_{0}^{s}\left [e^{(s-t)\Delta_{\mathbb{H}}}|\cdot|_{\mathbb{H}}^{-\frac{\gamma}{p-1}} \right]^{1-p}  dt \right\|_{L^\infty(\mathbb{H}^N)}  \\
\noalign{\medskip}
\geq  Cs^{-\frac{Q}{2}(p-1)} 
\left(\int_{B_{\mathbb{H}}(\eta_0,r)} e^{-{C_{\star} |\xi|_{\mathbb{H}}^2   }} u_0(\xi) d\xi \right)^{p-1} \cdot \\  
\hskip20pt \int_{0}^{s} (s-t)^{\frac{\gamma}{2}} \left[   \int_{\mathbb{H}^N} e^{-c_\star  |\eta |_{\mathbb{H}}^2   } |\eta |^{ - \frac{\gamma}{p-1}}_{\mathbb{H}}  d\eta   \right]^{1-p}  dt  \\
\noalign{\medskip}
\geq  Cs^{-\frac{Q}{2}(p-1)}  \int_{0}^{s} (s-t)^{\frac{\gamma}{2}} ds,
\end{array}
\end{equation*}
which contradicts condition \eqref{eq-22t11-b-neg}, letting $s\to +\infty$, since $1<p<1+(2+\gamma)/Q.$ 

The case $p=1+(2+\gamma)/Q$ can be treated as in the proof of Theorem \ref{global-e-blowup-gamapositivo}. \\

\textbf{ (ii)} Using the same argument as \cite[Lemma 2.12]{LEENI}, the estimate \eqref{estimativadokernel} and Lemma \ref{Lemachave1}, it is possible to verify that $\|e^{\sigma \Delta_{\mathbb{H}}}w_0\|_{L^\infty(\mathbb{H}^N)} \leq C {\sigma}^{-\frac{Q}{2}}$ for $\sigma>0$ and $w_0(\eta) = \phi(x) \phi(y) \theta(\tau) $ with $\phi (x) \leq C_1(1+|x|^2)^{-\frac{a}{2}}$, $a>N$ and $\theta(\tau) \leq C_2 (1+|\tau |^2)^{-\frac{1}{2}}$. Note that  if $0<t\leq 1$, then
$$\int_{0}^{t}  \|  e^{\sigma \Delta_{\mathbb{H}} }w_{0}  \|_{L^\infty(\mathbb{H}^N)}^{\frac{p-1}{q}} (t-\sigma)^{\frac{\gamma}{2}} d\sigma  \leq C \int_{0}^{t} (t-\sigma)^{\frac{\gamma}{2}} d\sigma \leq C t^{1+\frac{\gamma}{2}} \leq C.$$
Now, if $t\geq 1$, let $0<t_0<1/2$. Since  $t - \sigma \geq t/2 \ \mbox{ for } \ 0 \leq \sigma \leq t/2$, we obtain
\begin{align*}
    \int_{0}^{t}  \|  e^{\sigma \Delta_{\mathbb{H}} }w_{0}  \|_{L^\infty(\mathbb{H}^N)}^{\frac{p-1}{q}} (t-\sigma)^{\frac{\gamma}{2}} d\sigma &\leq Ct^{\gamma/2}+ C\left (\int_{t_0}^{t/2} + \int_{t/2}^t \right )  {\sigma}^{-\frac{Q(p-1)}{2q}} (t - \sigma)^{\frac{\gamma}{2}} d\sigma\\
    \noalign{\medskip}
    & \leq C+ 2^{-\gamma}Ct^{\frac{\gamma}{2}} \int_{t_0}^{t/2} {\sigma}^{-\frac{Q(p-1)}{2q}} d\sigma + Ct^{-\frac{Q(p-1)}{2q} + \frac{\gamma}{2} +1}\\
    \noalign{\medskip}
    & \leq C+Ct^{\frac{\gamma}{2}}  + Ct^{-\frac{Q(p-1)}{2q}+ \frac{\gamma}{2} +1}.
\end{align*}
Therefore, $\Lambda$ in \eqref{usaas} is verified if $-{Q(p-1)}/{2q} + \gamma/2 +1<0$, or $p\geq 1+ q {(2+\gamma)}/{Q}$. As $q>Q/(Q+\gamma)$, then $p\geq 1+ {(2+\gamma)}/{(Q+\gamma)}.$

~~\\
~~\\
\textbf{Acknowledgments}\\
~~\\
R. Castillo is supported by ANID-FONDECYT project No. 11220152. R. Freire is supported by CNPq (Conselho Nacional de Desenvolvimento Científico e Tecnológico) Process 140600/2020-5. M. Loayza is partially supported by MATH-AMSUD project 21-MATH-03 (CTMicrAAPDEs), CAPES-PRINT 88887.311962/2018-00 (Brazil), CNPq - 313382/2023-9.


\begin{thebibliography}{99}

\bibitem{Adams} R. A. Adams and J. J. Fournier, Sobolev spaces. Elsevier, 2003. 

\bibitem{AHMAD} B. Ahmad, A. Alsaedi, M. Kirane. Nonexistence of global solutions of some nonlinear space-nonlocal evolution equations on the Heisenberg group. Electron. J. Differential Equations, v. 10, (2015).

\bibitem{BandleLevine} C. Bandle,  H. A. Levine. On the existence and nonexistence of global solutions of reaction-diffusion equations in sectorial domains. Trans. Amer. Math. Soc. 316.2 (1989) 595-622.






\bibitem{BIRINDELLI-A} I. Birindelli, I. Dolcetta, A. Cutri. Liouville theorems for semilinear equations on the Heisenberg group. In: Annales de l'Institut Henri Poincaré C, Analyse non linéaire. No longer published by Elsevier, (1997) 295-308.


\bibitem{BONFIGLIOLI} A. Bonfiglioli, E. Lanconelli, F. Uguzzoni. Stratified Lie groups and potential theory for their sub-Laplacians. Springer Science \& Business Media, 2007.

\bibitem{breziscaz} H. Brezis, T. Cazenave. Nonlinear evolution equations. IM-UFRJ, Rio, v. 1, p. 994, 1994.

\bibitem{Ovidiu} O. Calin, D. Chang. Sub-Riemannian geometry: general theory and examples. Cambridge University Press, 2009.





\bibitem{CASTILLO1} R. Castillo, O. Guzmá-Rea, M. Loayza. On the local existence for Hardy parabolic equations with singular initial data. J. Math. Anal. Appl., v. 510, n. 2 (2022) 
126022. 

\bibitem{CastilloHardy} R. Castillo, R. Freire, M. Loayza, Global existence versus blow-up for a Hardy-Hénon parabolic equation on arbitrary domains. J. Differential Equations 429 (2025), 427-459.

\bibitem{CHI} N. Chikami, M. Ikeda, K. Taniguchi. Optimal well-posedness and forward self-similar solution for the Hardy–Hénon parabolic equation in critical weighted Lebesgue spaces. Nonlinear Analysis, v. 222,  (2022) 112931.

\bibitem{DAVIES1} E. Davies. Heat kernels and spectral theory. Cambridge university press, 1989.




\bibitem{FOLLAND1} G. B. Folland. Subelliptic estimates and function spaces on nilpotent Lie groups. Ark. Mat., v. 13, n. 1, (1975) 161-207.


\bibitem{Fujishima1} Y. Fujishima, N. Ioku,  Existence and nonexistence of solutions for the heat equation with a superlinear source term. J. Math. Pures Appl. (9) 118 (2018) 128-158.

\bibitem{Fujishima2} Y. Fujishima, K. Hisa, K. Ishige, R. Laister, Local solvability and dilation-critical singularities of supercritical fractional heat equations. J. Math. Pures Appl. (9) 186 (2024) 150-175.

\bibitem{FU} H. Fujita. On the blowing up of solutions of the Cauchy problem for $u_t = \Delta u+ u^{1+ \sigma}$. J. Fac. Sci. Univ. Tokyo, v. 13, (1966) 109-124.

\bibitem{FU1} H. Fujita, On some nonexistence and nonuniqueness theorems for nonlinear parabolic equations, Nonlinear Funct. Anal. (1970) 105-113.

\bibitem{PALMIERI} V. Georgiev, A. Palmieri. Lifespan estimates for local in time solutions to the semilinear heat equation on the Heisenberg group. Ann. Mat. Pura Appl., v. 200, n. 3 (2021) 999-1032.

\bibitem{GEORGIEV2} V. Georgiev, A. Palmieri. Critical exponent of Fujita-type for the semilinear damped wave equation on the Heisenberg group with power nonlinearity. J. Differential Equations, v. 269, n. 1, (2020) 420-448.

\bibitem{GOLDSTEIN} J. Goldstein, Q. Zhang. On a degenerate heat equation with a singular potential. J. Funct. Anal., v. 186, n. 2, p. (2001) 342-359.

\bibitem{Goldstein1}  G. R. Goldstein, J. Goldstein, A. E. Kogoj; A. Rhandi,  C. Tacelli. Instantaneous blowup and singular potentials on Heisenberg groups. Ann. Sc. Norm. Super. Pisa Cl. Sci. (5) 23 no. 4, (2022) 1723-1748.

\bibitem{Ghou} N. Ghoussoub and A. Moradifam, Functional inequalities: new perspectives and new applications, Mathematical Surveys and Monographs, vol. 187, American Mathematical Society, Providence, RI, 2013.

\bibitem{inequalities} G. H. Hardy, J. E. Littlewood, G. Pólya. Inequalities. Cambridge university press, 1952.

\bibitem{Henon} M. H\'enon, Numerical experiments on the stability of spherical stellar systems, Astron. Astrophys 24 (1973), 229–238.

\bibitem{Lars} L. H\"ormander. Hypoelliptic second order differential equations. 1967.

\bibitem{HULANICKI1} A. Hulanicki. Subalgebra of $L_1 (G)$ associated with laplacian on a Lie group, Colloquium Mathematicum (1974) 259-287.

\bibitem{Khare} K. Khare, M. Butola, S. Rajora, Fourier Optics and Computational Imaging; Springer Nature: Berlin/Heidelberg, Germany, 2023.

\bibitem{LEENI} T. Lee, W. Ni. Global existence, large time behavior and life span of solutions of a semilinear parabolic Cauchy problem. Trans. Amer. Math. Soc., v. 333, n. 1 (1992) 365-378.

\bibitem{Lieb1} E. Lieb, M. Loss, Analysis, Graduate Studies in Mathematics, American Mathematical Society, 2001.


\bibitem{Pascucci1} A. Pascucci. Semilinear equations on nilpotent Lie groups: global existence and blow-up of solutions. Matematiche (Catania) 53, no. 2 (1998) 345-357.


\bibitem{PINSKY1} R. Pisnky. Existence and Nonexistence of Global Solutions for $u_t =  \Delta u +a(x)u^p \in \mathbb{R}^d$. J. Differential Equations, v. 133, n. 1 (1997) 152-177.

\bibitem{MITIDIERI} S. Pohozaev, L. Véron. Nonexistence results of solutions of semilinear differential inequalities on the Heisenberg group. Manuscripta Math., v. 102, n. 1 (2000).

\bibitem{QI} Y. Qi. The critical exponents of parabolic equations and blow-up in Rn. Proc. Roy. Soc. Edinburgh Sect. A, v. 128, n. 1, (1998) 123-136.

\bibitem{QS} P. Quittner, P. Souplet. Superlinear parabolic problems. Springer International Publishing, 2019.

\bibitem{RANDALL} J. Randall. The heat kernel for generalized Heisenberg groups. J. Geom. Anal., v. 6, (1996) 287-316.

\bibitem{RUZHANSKY} M. Ruzhansky, N. Yessirkenov. Existence and non-existence of global solutions for semilinear heat equations and inequalities on sub-Riemannian manifolds, and Fujita exponent on unimodular Lie groups. J. of Differential Equations, v. 308 (2022) 455-473.

\bibitem{SL} B. Slimene, S. Tayachi, F. Weissler. Well-posedness, global existence and large time behavior for Hardy–Hénon parabolic equations. Nonlinear Analysis, v. 152, (2017) 116-148.




\bibitem{VAROPOULOS} N. T. Varapoulos, L. Sallof-Coste, T. Coulhon. Cambridge Tracts in Math. Analysis and Geometry on Groups, v. 100, 1993.

\bibitem{WANG} X. Wang. On the Cauchy problem for reaction-diffusion equations. Trans. Amer. Math. Soc., v. 337, n. 2, (1993) 549-590.


\bibitem{Weissler} F. Weissler. Existence and nonexistence of global solutions for a semilinear heat equation. Israel J. Math. 38 no. 1-2 (1981) 29-40.



\bibitem{Weissler1} F. Weissler. Local existence and nonexistence for semilinear parabolic equations in $L\sp{p}$. Indiana Univ. Math. J. 29 no. 1, (1980) 79-102.


 
\bibitem{Weissler2} F. Weissler. $L^p$-energy and blow-up for a semilinear heat equation. Nonlinear functional analysis and its applications,  Proc. Sympos. Pure Math., 45, Part 2, Amer. Math. Soc., Providence, RI (1986) 545-551.

    

\end{thebibliography}
\end{document}